\documentclass[12pt]{amsart}
\usepackage{amscd,amsmath,amssymb,amsfonts}
\usepackage[cmtip, all]{xy}

\newtheorem{thm}{Theorem}[section]
\newtheorem{prop}[thm]{Proposition}
\newtheorem{lem}[thm]{Lemma}
\newtheorem{cor}[thm]{Corollary}
\newtheorem{conj}[thm]{Conjecture}

\theoremstyle{definition}
\newtheorem{exm}[thm]{Example}
\newtheorem{defn}[thm]{Definition}

\theoremstyle{remark}
\newtheorem{remk}[thm]{Remark}
\newtheorem{remks}[thm]{Remarks}
\newtheorem{exms}[thm]{Examples}
\newtheorem{notat}[thm]{Notation}

\numberwithin{equation}{section}

{\hfill$\square$\end{defn}}

{\hfill$\square$\end{remk}}

{\hfill$\square$\end{remks}}

{\hfill$\square$\end{exm}}

{\hfill$\square$\end{exms}}

{\hfill$\square$\end{notat}}

\newcommand{\sC}{{\mathcal C}}

\newcommand{\sH}{{\mathcal H}}
\newcommand{\sI}{{\mathcal I}}
\newcommand{\sJ}{{\mathcal J}}
\newcommand{\sK}{{\mathcal K}}
\newcommand{\sL}{{\mathcal L}}

\newcommand{\sO}{{\mathcal O}}

\newcommand{\wt}{\widetilde}

\newcommand{\A}{{\mathbb A}}

\newcommand{\C}{{\mathbb C}}

\renewcommand{\P}{{\mathbb P}}
\newcommand{\Q}{{\mathbb Q}}

\newcommand{\Z}{{\mathbb Z}}

\newcommand{\surj}{\twoheadrightarrow}
\newcommand{\inj}{\hookrightarrow}

\newcommand{\del}{\partial}

\newcommand{\ds}{{/\kern-3pt/}}

\newcommand{\ov}{\overline}

\begin{document}
\title{An Artin-Rees Theorem and
Applications to Zero Cycles}
%\title{An Artin-Rees Theorem in $K$-theory and
%Applications to Zero Cycles}
\author{Amalendu Krishna}


\baselineskip=10pt 
                                              
\begin{abstract}
For the smooth normalization $f : {\ov X} \to X$ of a singular variety 
$X$ over a field $k$ of characteristic zero, we show that for any 
conducting subscheme $Y$ for the normalization, and for any $i \in \Z$,
the natural map $K_i(X, {\ov X}, nY) \to K_i(X, {\ov X}, Y)$ is zero
for all sufficiently large $n$.

As an application, we prove a formula for the Chow group of zero cycles
on a quasi-projective variety $X$ over $k$ with Cohen-Macaulay isolated
singularities, in terms of an inverse limit of relative Chow groups of a 
desingularization $\wt X$ relative to the multiples of the exceptional
divisor.

We use this formula to verify a conjecture of Srinivas about the
Chow group of zero cycles on the affine cone over a smooth projective
variety which is arithmetically Cohen-Macaulay.
\end{abstract} 

\maketitle

\section{Introduction}
It has been shown over the years that the relative and the double relative
$K$-groups play a very important role in the study of algebraic $K$-theory 
and algebraic cycles on singular varieties. These are particularly very
useful tools in comparing the $K$-groups (or Chow groups) of a singular
scheme and its normalization. The double relative $K$-groups are in
general very difficult to compute. Our aim in this paper is two fold.
First we prove an Artin-Rees type theorem for the double relative $K$-groups,
and then give some very important applications of this result to convince
the reader why such Artin-Rees type results are indeed desirable if one
wants to study algebraic cycles on singular varieties.

In this paper, a {\sl variety} will mean a reduced, connected and
separated scheme of finite type over a field $k$ of characteristic
zero. All the underlying fields in this paper will be of characteristic zero.
Let $X$ be singular variety of dimension $d$ over a field $k$ and let
$f : {\ov X} \to X$ be the normalization morphism. For any closed subscheme 
$Y$ of $X$, the relative $K$-groups $K_*(X, Y)$ are defined as the stable 
homotopy groups of the homotopy fiber $K(X, Y)$ of the map $K(X) \to K(Y)$ of 
non-connective spectra as defined in \cite{TT}.
If $Y \inj X$ is a conducting subscheme for the normalization, we put
${\ov Y} = Y {\times}_X {\ov X}$. We define the double relative $K$-groups 
$K_*(X, {\ov X}, Y)$ as the stable homotopy groups of the homotopy fiber 
$K(X, {\ov X}, Y)$ of the induced natural map $K(X, Y) \to 
K({\ov X},{\ov Y})$ of non-connective spectra. As shown in \cite{TT},
the underlying spectra to define the relative and double relative $K$-groups
are non-connective, and hence the groups $K_i(X, {\ov X}, Y)$ may be
non-zero even when $i$ is negative. For a conducting subscheme $Y$, we
let $\sI$ be the sheaf of ideals on $X$ defining $Y$. 
We denote by $nY$ (for $n \ge 1)$, the conducting subscheme of $X$ for the 
normalization that is defined by the sheaf of ideals ${\sI}^n$. The 
functoriality of the double relative groups defines the natural maps 
$K_*(X, {\ov X}, (n+1)Y) \to K_*(X, {\ov X}, nY)$ for all $n \ge 1$.

If $X$ is affine, it is classically known that $K_i(X, {\ov X}, Y) = 0$
for $i \le 0$. Furthermore, it follows from \cite{CGW} (Theorem~3.6) and a 
result of Cortinas (\cite{Cortinas}, Corollary~0.2) that $K_1(X, {\ov X}, Y)
\cong {{\sI}_Y}/{{\sI}^2_Y} \otimes {\Omega}_{{\ov Y}/Y}$. This shows that
for $X$ affine, the map $K_i(X, {\ov X}, 2Y) \to K_i(X, {\ov X}, Y)$ is zero 
for $i \le 1$. This has already been very useful in the study of zero cycles
on singular schemes. We refer the reader to \cite{KS}, \cite{Kr1}, 
\cite{BPW} among several others for some applications of such results to 
zero cycles on surfaces and threefolds.  
One of the main difficulties in advancing the study of $K$-theory and
algebraic cycles to the higher dimensional singular varieties has been the 
need for the following general Artin-Rees type result.   
\begin{thm}\label{thm:Artin-Rees}
Let $X$ be a quasi-projective variety of dimension $d$ over a field $k$
and let $f: {\ov X} \to X$ be the normalization of $X$. Assume that 
$\ov X$ is smooth. Then for any conducting subscheme $Y \inj X$ for the 
normalization, and for any $i \in \Z$, the natural map $ K_i(X, {\ov X}, nY) 
\to K_i(X, {\ov X}, Y)$ is zero for all sufficiently large $n$.
\end{thm}

As already said before, one of the motivations for proving the above theorem
is its applications in computing the $K$-groups and the Chow groups of
algebraic cycles on singular varieties. These Chow groups are usually very
difficult to compute, even for varieties with isolated singularities. 
In this paper, we use the above
theorem to prove the following formula for the Chow group of zero cycles on 
a quasi-projective variety $X$ with Cohen-Macaulay isolated singularities,
in terms of an inverse limit of the relative Chow groups of a 
desingularization $\wt X$ relative to the multiples of the exceptional
divisor. Such a formula for the normal surface singularities was 
conjectured by Bloch and Srinivas (\cite{Srinivas1}) and proved by 
Srinivas and the author in \cite{KS} (Theorem~1.1). This formula was
later verified for the threefolds with Cohen-Macaulay isolated singularities
in \cite{Kr1} (Theorem~1.1). 

Before we state our next theorem, we very briefly recall the definition
of the Chow group of zero cycles on singular varieties and some other
related notions. Let $X$ be a quasi-projective singular variety of
dimension $d$ over a field $k$. The cohomological 
{\sl Chow group of zero cycles} $CH^d(X)$ was defined by Levine and
Weibel in {\cite{LevineW} as the free abelian group on the smooth closed
points of $X$, modulo the zero cycles given as the sum of divisors of suitable
rational functions on certain {\sl Cartier curves} on $X$. These curves
can be assumed to be irreducible and disjoint from the singular locus of
$X$, if it is normal. We refer the reader to \cite{ESV} for the full
definition of $CH^d(X)$ and some of its important properties. 

Let $F^{d}K_0(X)$ denote the subgroup of the Grothendieck group of vector 
bundles $K_0(X)$, generated by the classes of smooth codimension $d$ closed 
points on $X$. There is a natural surjective map 
$CH^d(X) \rightarrow F^dK_0(X)$. This map is known ({\sl cf.} \cite{Levine1}, 
Corollary~2.6 and Theorem~3.2) to be an isomorphism, if the underlying field 
$k$ is algebraically closed and $X$ is either an affine or a projective 
variety over $k$ with only normal isolated singularities. This isomorphism 
will be used throughout this paper.  
For any closed subscheme $Z$ of $X$, let $F^{d}K_{0}(X,Z)$ be the subgroup of 
$K_{0}(X,Z)$ (defined above) generated by the classes of 
smooth closed points of $X-Z$ ({\sl cf.} \cite{CoombesS}). We shall also call
$F^{d}K_{0}(X,Z)$ the {\sl relative Chow group of zero cycles} of $X$ 
relative to $Z$. 

Now let $X$ be an irreducible quasi-projective variety of dimension $d$ with 
isolated singularities over a field $k$, and let $p : {\wt X} \to X$ be a 
resolution of singularities of $X$.
Let $E$ denote the reduced exceptional divisor on $\wt X$ and
let $nE$ denote the $n$th infinitesimal thickening of $E$. From the
above definitions and from Lemma~\ref{lem:commute} below
(see also \cite{Kr1}, Theorem~1.1), one has the following commutative 
diagram for each $n > 1$ with all arrows surjective.
\begin{equation}\label{eqn:relative}
\xymatrix{
&& F^{d}K_{0}(X) \ar @{->>} [dll] \ar @{->>} [d] \ar @{->>}[drr]^{p^{*}}
&& \\ F^{d}K_{0}({\wt X}, nE) \ar @{->>} [rr] &&
F^{d}K_{0}({\wt X} , (n-1)E)
\ar @{->>} [rr] && F^{d}K_{0}({\wt X}) 
}
\end{equation}

We also recall from \cite{Kr1} that a resolution of singularities
$p : {\wt X} \to X$ of $X$ is called a {\sl good resolution
of singularities} if $p$ is obtained as a blow up of
$X$ along a closed subscheme $Z$ whose support is the set of singular
points of $X$. Such a resolution of singularities always exists.
In fact, it is known ({\sl cf.} \cite{Kr1}, Lemma~2.5) using Hironaka's 
theory that every resolution of singularities is dominated by a good
resolution of singularities. Hence in the study of the Chow group of zero
cycles on $X$, it suffices to consider only the good resolutions
of singularities since $CH^d(\wt X)$ is a birational invariant of 
a smooth projective variety $\wt X$ of dimension $d$. We now state our next 
result.
\begin{thm}\label{thm:Chowgroup}
Let $X$ be a quasi-projective variety of dimension $d \ge 2$ over a field $k$. 
Let $p : {\wt X} \to X$ be a good resolution of singularities of $X$
with the reduced exceptional divisor $E$ on $\wt X$. 
Assume that $X$ has only Cohen-Macaulay isolated singularities.
Then for all sufficiently large $n$, the maps \\
$(i)$ $F^{d}K_{0}(\wt {X},nE) \to F^{d}K_{0}(\wt {X},(n-1)E)$ and \\
$(ii)$ $F^{d}K_{0}(X) \to F^{d}K_{0}(\wt {X},nE)$ \\
are isomorphisms. In particular, if $k = {\ov k}$ and $X$ is either affine
or projective, then 
\[CH^d(X) \cong~ ~{\underset {n}{\varprojlim}}~F^{d}K_{0}(\wt {X},nE).\]
\end{thm}	
Since rational singularities are Cohen-Macaulay, we immediately get
\begin{cor}\label{cor:rational}
If $X$ is a quasi-projective variety with only rational and isolated 
singularities, then one has the above formula for $CH^d(X)$.
\end{cor}  
It turns out that the relative Chow groups (as defined above) of the 
resolution of singularities, relative to the multiples
of the exceptional divisor are often computable. Such computations have
been carried out before to compute the Chow groups of zero cycles on normal
surfaces ({\sl cf.} \cite{KS}) and certain threefolds ({\sl cf.} \cite{Kr1}).
We now proceed to carry these computations further in all dimensions and 
then use Theorem~\ref{thm:Chowgroup} to compute the Chow group of zero 
cycles on certain classes of singular varieties. These results in
particular give many new examples of the validity of the singular analogue 
of the well known conjecture Bloch and its generalization in higher
dimensions. 

Our result in this direction is inspired by the following conjecture of 
Srinivas (\cite{Srinivas2}, Section~3), which in itself was probably
motivated by the affine analogue of a conjecture of Bloch about the Chow 
groups of zero cycles on smooth projective varieties.
\begin{conj}[Srinivas]\label{conj:Srinivas}
Let $Y \inj {\P}^N_{\C} $ be a smooth and projectively normal variety over 
$\C$ of dimension $d$ and let $X = C(Y)$ be the affine cone over $X$. Then 
$CH^{d+1}(X) = 0$ if and only if $H^d(Y, {{\sO}_Y}(1)) = 0$.
\end{conj}
The `only if' part of this conjecture was proved by Srinivas himself in
\cite{Srinivas2} (Corollary~2). So the `if' part remains open at present. 
This conjecture was verified in \cite{KS} (Corollary~1.4) when $Y$ is
a curve. This was also verified by Consani (\cite{Consani}, Proposition~3.1)
in a very special case when $Y$ is a hypersurface in ${\P}^3$.
The following result verifies this conjecture when the affine cone $X$ is 
Cohen-Macaulay. In particular, Srinivas' conjecture is verified when
$Y$ is a complete intersection in ${\P}^N$.
\begin{thm}\label{thm:cone}
Let $k$ be an algebraically closed field of characteristic zero and let 
$Y \inj {\P}^N_{k}$ be a smooth projective variety of dimension $d$. Let 
$X = C(Y)$ be the affine cone over $Y$. Assume that $X$ is Cohen-Macaulay. 
Then $CH^{d+1}(X) = 0$ if $H^d(Y, {{\sO}_Y}(1)) = 0$. 

Moreover, if 
$k$ is a universal domain and if ${\ov X}$ is the projective cone over $Y$,
then the following are equivalent. \\
$(i)$ $H^d(Y, {{\sO}_Y}(1)) = 0$ \\
$(ii)$ $H^{d+1}({\ov X}, {\sO}_{\ov X}) = 0$ \\
$(iii)$ $H^{d+1}({\ov X}, {\Omega}^i_{{\ov X}/k}) \xrightarrow {\cong}
H^{d+1}(W, {\Omega}^i_{{W}/k})$ for all $i \ge 0$, where $W$ is a resolution 
of singularities of $\ov X$. \\
$(iv)$ $CH^{d+1}(X) = 0$ \\
$(v)$ $CH^{d+1}(\ov X) \cong CH^{d}(Y)$.
\end{thm}
If $Y \inj {\P}^N_{k}$ is a hypersurface of degree $d$, then there is an
isomorphism \\
$H^{N-1}\left(Y, {{\sO}_Y}(1)\right) \cong 
H^0\left({\P}^N_k, {\sO}_{{\P}^N} (d-N-2)\right)$ and we obtain 
the following immediate consequence of Theorem~\ref{thm:cone}.
\begin{cor}\label{cor:rational1}
Let $Y \inj {\P}^N_{k}$ be a smooth hypersurface of degree $d$ and let
$X = C(Y)$ be the affine cone over $Y$. Then $CH^{N}(X) = 0$ if $d \le N+1$.
The converse also holds if $k$ is a universal domain.
\end{cor}  
The following application of Theorem~\ref{thm:cone} to the projective modules
on singular affine algebras follows at once from
Murthy's result (\cite{Murthy}, Corollary~3.9).
\begin{cor}\label{cor:projective}
Let $Y  \inj {\P}^N_{k}$ be as in Theorem~\ref{thm:cone} and let $A$ be 
its homogeneous coordinate ring. If $H^d(Y, {{\sO}_Y}(1)) = 0$, then
every projective module over $A$ of rank at least $d$ has a unimodular element.
\end{cor}  
If $X$ is a normal projective surface with a resolution of singularities
$p : {\wt X} \to X$, then it was shown in \cite{KS} (Theorem~1.3) that
$CH^2(X) \cong CH^2(\wt X)$ iff $H^2(X, {\sO}_X) \cong 
H^2({\wt X}, {\sO}_{\wt X})$. If $X$ has higher dimension, it is not expected
that the isomorphism of the top cohomology of the structure sheaves
is the sufficient condition to conclude the isomorphism of the Chow groups
of zero cycles. However, Theorem~\ref{thm:cone} suggests that this 
might still be true in the case of isolated singularities.
%\begin{conj}\label{conj:isolated}
%Let $X$ be a normal projective variety of dimension $d$ with isolated 
%%singularities and let $p: {\wt X} \to X$ be a resolution of singularities of 
%$X$. Then $CH^d(X) \cong CH^d(\wt X)$ if $H^{d}(X, {\sO}_X) \cong
%H^{d}({\wt X}, {\sO}_{\wt X})$.
%\end{conj}
%As our final application of Theorem~\ref{thm:Chowgroup}, we obtain the 
%following result about the Chow group of zero cycles on varieties with
%quotient singularities.
%\begin{cor}\label{cor:quotient}
%Let $G$ be a finite group acting on a smooth projective variety $V$ of
%dimension $d$ over $k$ and let $X = V/G$ be the quotient. Let
%$p: {\wt X} \to X$ be a resolution of singularities of $X$. Then
%the map $CH^d(X) \to CH^d(\wt X)$ is an isomorphism.
%\end{cor} 

As can be seen in the above results, our applications of 
Theorem~\ref{thm:Artin-Rees} has been restricted to computing some
pieces of the group $K_0$ of singular varieties. On the other hand, the
conclusion of this theorem applies also to higher $K_i$'s. We hope that
this result will be a significant tool to study $K_1$ of singular 
varieties for which hardly anything is known.

We conclude this section with a brief outline of this paper.
Our strategy of proving Theorem~\ref{thm:Artin-Rees} is to use
the Brown-Gersten spectral sequence of \cite{TT} and the Cortinas'
proof of KABI-conjecture in \cite{Cortinas} to reduce the problem
to proving similar results for the Hochschild and cyclic homology.
Since these homology groups have canonical decompositions in terms
of Andr{\'e}-Quillen homology, we first prove an Artin-Rees theorem
for these homology groups. We give an overview of the Andr{\'e}-Quillen 
homology and then generalize a result of Quillen (\cite{Quillen}, Theorem~6.15)
in the next section. Section~3 contains a proof of our Artin-Rees type
result for the Andr{\'e}-Quillen homology. In Section~4, we deal with
proving some refinements of this result for the special case of
conducting ideals for the smooth normalization of the essentially of finite
type $k$-algebras, where $k$ is a field. Section~5 generalizes these
results for the Andr{\'e}-Quillen homology over any base field which
is contained in $k$. We then use these results in Section~6 to prove the 
Artin-Rees theorem for the double relative Hochschild and cyclic homology,
and then give the proof of Theorem~\ref{thm:Artin-Rees}. Section~7 contains
the proof of Theorem~\ref{thm:Chowgroup}. In Section~8, we relate
the Chow group of zero cycles with certain cohomology of Milnor $K$-sheaves
and then prove some general results about these cohomology groups, which
are then used in the final section to compute the Chow group of zero cycles 
on the affine cones. 

\section{Andr{\'e}-Quillen Homology}
All the rings in this section will be assumed to be commutative $k$-algebras, 
where $k$ is a given field of characteristic zero. Our aim in this section is
to give an overview of Andr{\'e}-Quillen homology and related concepts.
This homology theory of algebras will be our main tool to prove 
Theorem~\ref{thm:Artin-Rees}. We also prove here 
a generalization of an Artin-Rees type theorem of Quillen
{\cite{Quillen}, Theorem~6.15) for the Andr{\'e}-Quillen homology of 
algebras which are finite over the base ring.
As Quillen shows (see also \cite{Andre}), such a result itself has many 
interesting consequences for the homology of commutative algebras.
Apart from the above cited works of
Andr{\'e} and Quillen, our other basic reference for this material
including the Hochschild and cyclic homology, is \cite{Loday}.

Let $A$ be a commutative ring which is essentially of finite type over
the field $k$, and let $B$ be an $A$-algebra. A simplicial $A$-algebra
will mean a simplicial object in the category of $A$-algebras. Let $P_*$
be a simplicial $A$-algebra. We say that $P_*$ is $B$-augmented if the 
natural map $A \to B$ factors through $A \to P_* \to B$, where any 
$A$-algebra $B$ is naturally considered a simplicial $A$-algebra with
all the face and degeneracy maps taken as identity map of $B$. The homotopy 
groups of a simplicial $A$-algebra $P_*$ is defined as the homotopy groups
of the simplicial set $P_*$, which is same as the homotopy groups of the
simplicial $A$-module $P_*$. The Dold-Kan correspondence implies that
these homotopy groups are same as the homology groups of the corresponding
chain complex (which we also denote by $P_*$) of $A$-modules. Using this
equivalence between simplicial $A$-modules and chain complexes of
$A$-modules, we shall often write the homotopy groups ${\pi}_i(M_*)$
of a simplicial $A$-module $M_*$ homologically as $H_i(M_*)$ without any ado.
We say that $P_*$ is a {\sl free} $A$-algebra if each $P_i$ ($i \ge 0$)
is a symmetric algebra over a free $A$-module. 
\begin{defn}\label{defn:resolution}
A free simplicial $A$-algebra $P_*$ is called a {\sl simplicial resolution}
of an $A$-algebra $B$ if $P_*$ is $B$-augmented such that the natural 
map $H_i(P_*) \to H_i(B)$ is an isomorphism for all $i$.
\end{defn} 
It is known ({\sl cf.} \cite{Loday}, Lemma~3.5.2) that any $A$-algebra $B$
admits a free simplicial resolution and any two such resolutions are
homotopy equivalent. 

Before we define the cotangent modules, we recall 
({\sl cf.} \cite{Loday}, 1.6.8) that for two simplicial $B$
modules $M_*$ and $N_*$, their tensor and wedge products are defined
degree-wise, i.e., 
\[(M_* {\otimes}_B N_*)_i = M_i {\otimes}_B N_i \ {\rm and} \
\left({{\wedge}^r_B} M_*\right)_i =  {{\wedge}^r_B} M_i \ {\rm for} \
r, i \ge 0.\]
The face and degeneracy maps of the tensor (or wedge) product are
degree-wise tensor (or wedge) product of the corresponding maps. 
Since we are in characteristic zero, the following lemma relating tensor and 
exterior powers of a simplicial $B$-module is elementary.
\begin{lem}\label{lem:elem}
For any simplicial $B$-module $M_*$ and for $r \ge 0$, ${{\wedge}^r_B} M_*$
is canonically a retract of ${{\otimes}^r_B} M_*$. In particular,
$H_i\left({{\wedge}^r_B} M_*\right)$ is a canonical direct summand
of $H_i\left({{\otimes}^r_B} M_*\right)$ for all $i \ge 0$. 
\end{lem}
\begin{proof} This is well known and we only give a very brief sketch.
Since the tensor and exterior powers are defined degree-wise,
it suffices to prove the lemma for a $B$-module $M$. One defines a
$B$-linear map \[{\otimes}^r M \xrightarrow {alt} {\otimes}^r M \]
as
\[alt(a_1 \otimes \cdots \otimes a_r) = 1/{r!} 
\underset {\sigma \in S_r}{\Sigma} sgn(\sigma) a_{\sigma (1)} \otimes
\cdots \otimes a_{\sigma (r)}.\]
It is easy to check that `$alt$' is a projector and is natural with respect
to map of $B$-modules. Moreover, if ${\wt {\wedge}}^r M$ denotes the
image of this map, then the composite ${\wt {\wedge}}^r M \to
{\otimes}^r M \to {{\wedge}^r} M$ is a canonical isomorphism.
\end{proof}

For an $A$-algebra $B$, we define the {\sl cotangent module} of $B$ to be the
simplicial $B$-module ${\mathbb L}_{B/A}$ given by 
\begin{equation}\label{eqn:cotangent}
\left({\mathbb L}_{B/A}\right)_i = {\Omega}^1_{{P_i}/A} {\otimes}_{P_i} B,
\end{equation}
where $P_*$ is any free simplicial resolution of $B$ and the face and 
degeneracy maps of ${\mathbb L}_{B/A}$ are induced by those of $P_*$.
The homotopy equivalence of different simplicial resolutions of $B$
implies that ${\mathbb L}_{B/A}$ is a well-defined object in the homotopy 
category of simplicial $B$-modules. The {\sl Andr{\'e}}-{\sl Quillen homology}
of $B$ with coefficients in a $B$-module $M$ is defined as 
\[D_q(B/A , M) = H_q\left({\mathbb L}_{B/A} {\otimes}_B M\right), q \ge 0.\]

One defines the {\sl higher Andr{\'e}}-{\sl Quillen homology} of $B$
with coefficients in a $B$-module $M$ as 
\[D^{(r)}_q(B/A, M) = H_q\left({\mathbb L}^{(r)}_{B/A} {\otimes}_B M\right) \ 
{\rm for} \ r, q \ge 0,\]
where \[{\mathbb L}^{(r)}_{B/A} = {\wedge}^r_B\left({\mathbb L}_{B/A}\right) \
{\rm for} \ r \ge 0. \]
When $M$ is same as $B$, we write $D^{(r)}_q(B/A, M)$ simply as
$D^{(r)}_q(B/A)$. It is easy to see from these definitions that
$D^{(0)}_0(B/A, M) = M$ and $D^{(0)}_{\ge 1}(B/A, M) = 0$. It is
also known ({\sl cf.} \cite{Loday}, Theorem~3.5.8, Theorem~4.5.12)
that for any $r \ge 0$, there is a canonical isomorphism
\begin{equation}\label{eqn:D0}
D^{(r)}_0(B/A, M) \xrightarrow {\cong} {\Omega}^r_{B/A} {\otimes}_B M.
\end{equation}

Let $A$ be a $k$-algebra as above and let ${\{B_n\}}_{n \ge 0}$ be an
inverse system of $A$-algebras. We denote this inverse system by 
$B_{\bullet}$. A $B_{\bullet}$-module is an inverse system 
${\{M^n\}}_{n \ge 0}$ of $A$-modules such that for each $n \ge 0$,
$M^n$ is in fact a $B_n$-module and the map $M^n \xrightarrow {f^n} M^{n-1}$
is $B_n$-linear such that these maps are compatible with maps in
the inverse system $B_{\bullet}$. A typical example in which we shall be
mostly interested in is when $I$ is an ideal of $A$ and $M^n = B_n = 
A/{I^{n+1}}$.
\begin{prop}\label{prop:tensor}
Let $B_{\bullet}$ be an inverse system of $A$-algebras. Let $M^{\bullet}_*$
and $N^{\bullet}_*$ be the flat simplicial $B_{\bullet}$-modules such that
for each $q \ge 0$ and for each $n_0 \ge 0$, the map 
$H_q(M^n_*) \to H_q(M^{n_0}_*)$ is zero for all $n \gg n_0$. Then for each 
$q \ge 0$ and for each $n_0 \ge 0$, the map 
$H_q\left(M^n_* {\otimes}_{B_n} N^n_*\right) \to
H_q\left(M^{n_0}_* {\otimes}_{B_{n_0}} N^{n_0}_*\right)$ is zero for
all $n \gg n_0$.
\end{prop}
\begin{proof} For any $A$ algebra $B$, let $C: SimpMod (B) \to Ch_{\ge 0}(B)$
be the Dold-Kan functor from the category of simplicial $B$-modules to the 
category of chain complexes of $B$-modules which are bounded below at zero. 
This functor takes a simplicial module $M_*$ to itself and the differential
at each level is the alternating sum of the face maps at that level.
Then the Eilenberg-Zilber theorem ({\sl cf.} \cite{Loday}, 1.6.12)
implies that there is a natural Alexander-Whitney map 
\[C\left(M_* {\otimes}_B N_*\right) \to C(M_*) {\otimes}_B C(N_*)\] which is a 
quasi-isomorphism. Here the term on the right is the tensor product in the
category of chain complexes, which is given as the total complex of the
double complex $(M_*, N_*)_{i, j} = M_i {\otimes}_B N_j$. Hence it suffices
to prove the lemma for the tensor product of chain complexes.

Since the double complex $\{M^n_i {\otimes}_{B_n} N^n_j\}_{i, j \ge 0}$ lies
only in the first quadrant, there is a convergent spectral sequence 
(\cite{Weibel1}, 5.6.1)  
\[{^n}E^2_{p,q} = H_p\left(\cdots \to H_q(M^n_* {\otimes}_{B_n} N^n_j) \to
\cdots \to H_q(M^n_* {\otimes}_{B_n} N^n_1) \to 
H_q(M^n_* {\otimes}_{B_n} N^n_0)\right)\] 
\[\hspace*{7cm}\Rightarrow H_{p+q}\left(M^n_* {\otimes}_{B_n} N^n_*\right).\]
This spectral sequence is compatible with the maps in the inverse systems
$\{B_n\}$, $\{M^n\}$ and $\{N^n\}$ and we get an inverse system of
spectral sequences
\[
\xymatrix{
{^n}E^2_{p,q} \ar[d] & {\Rightarrow} & 
H_{p+q}\left(M^n_* {\otimes}_{B_n} N^n_*\right) \ar[d] \\
{^{n-1}}E^2_{p,q} & {\Rightarrow} & 
H_{p+q}\left(M^{n-1}_* {\otimes}_{B_{n-1}} N^{n-1}_*\right).}
\]
Now as $N^n_j$ is a flat $B_n$-module for each $n$ and $j$, we see that
${^n}E^2_{p,q}$ is same as 
$H_p\left(H_q(M^n_*) {\otimes}_{B_n} N^n_*\right)$. 
In particular, we see that for each $p, q \ge 0$ and for each $n_0 \ge 0$,
the natural map
${^n}E^2_{p,q} \to {^{n_0}}E^2_{p,q}$ is zero for all $n \gg n_0$ and
hence the map ${^n}E^i_{p,q} \to {^{n_0}}E^i_{p,q}$ is zero for all 
$n \gg n_0$ and for all $i \ge 2$. From this we conclude that for a fixed
$q \ge 0$, there is a filtration 
\[0 = F^n_{-1} \subset F^n_0 \subset \cdots F^n_{q-1} \subset 
F^n_q \left(H_q \left(M^n_* {\otimes}_{B_n} N^n_*\right)\right) =
H_q \left(M^n_* {\otimes}_{B_n} N^n_*\right)\]
and a map of filtered $B_n$-modules 
$H_q \left(M^n_* {\otimes}_{B_n} N^n_*\right) \to
H_q \left(M^{n-1}_* {\otimes}_{B_{n-1}} N^{n-1}_*\right)$ such that
for each $j, n_0 \ge 0$, the map ${F^n_j}/{F^n_{j-1}} = {^n}E^{\infty}_{j, q-j}
\to {^{n_0}}E^{\infty}_{j, q-j} = {F^{n_0}_j}/{F^{n_0}_{j-1}}$ is zero for all
$n \gg n_0$. 

Now we show by induction that for any $0 \le j \le q$ 
and any $n_0 \ge 0$, the map $F^n_j \to F^{n_0}_j$ is zero for all
$n \gg n_0$. This will finish the proof of the proposition. The following 
trick (which we call the {\sl doubling trick}) to do this will be used 
repeatedly in this paper. We fix $j$ with 
$0 \le j \le q$ and by induction, assume that
there exist $n_1 \gg n_0$ and $n_2 \gg n_1$ such that in the commutative
diagram 
\begin{equation}\label{eqn:doubling}
\xymatrix{
0 \ar[r] & F^n_{j-1} \ar[d] \ar[r] & F^n_j \ar[r] \ar[d] &
{^n}E^{\infty}_{j, q-j} \ar[d] \ar[r] & 0 \\
0 \ar[r] & F^{n_1}_{j-1} \ar[d] \ar[r] & F^{n_1}_j \ar[r] \ar[d] &
{^{n_1}}E^{\infty}_{j, q-j} \ar[d] \ar[r] & 0 \\
0 \ar[r] & F^{n_0}_{j-1} \ar[r] & F^{n_0}_j \ar[r] &
{^{n_0}}E^{\infty}_{j, q-j} \ar[r] & 0, }
\end{equation}
the bottom left and the bottom right vertical arrows are zero for all
$n \ge n_1$, and the top left and the top right vertical arrows are zero
for all $n \ge n_2$. A diagram chase shows that the composite middle
vertical arrow is zero for all $n \ge n_2$.
\end{proof}

\begin{cor}\label{cor:exterior}
Let $B_{\bullet}$ be an inverse system of $A$-algebras and let 
$M^{\bullet}_*$ be a flat simplicial $B_{\bullet}$-module. Assume that for 
each $q \ge 0$ and for each $n_0 \ge 0$, the map $H_q(M^n_*) \to 
H_q(M^{n_0}_*)$ is zero for all $n \gg n_0$. Then for each $r \ge 1$ and 
$q , n_0 \ge 0$, the map
$H_q\left({\wedge}^r_{B_n} M^n_*\right) \to  
H_q\left({\wedge}^r_{B_{n_0}} M^{n_0}_*\right)$ is zero for all $n \gg n_0$.
\end{cor}
\begin{proof} By Lemma~\ref{lem:elem}, it suffices to prove the corollary
when the exterior powers are replaced by the corresponding tensor powers, 
when it follows directly from Proposition~\ref{prop:tensor} and induction
on $r$.
\end{proof}
The following result was proved by Quillen (\cite{Quillen}, Theorem~6.15)
for $r = 1$.
\begin{cor}\label{cor:AQ1}
Let $A$ be an essentially of finite type $k$-algebra. Let $I$ be an ideal
of $A$ and put $B_n = A/{I^{n+1}}$ for $n \ge 0$. Then for each 
$r \ge 1$ and  $q, n_0 \ge 0$, the natural map $D^{(r)}_q({B_n}/A) \to
D^{(r)}_q({B_{n_0}}/A)$ is zero for all $n \gg n_0$.
\end{cor}  
\begin{proof} We see from ~\ref{eqn:cotangent} that ${\mathbb L}_{{B_n}/A}$ 
is a free simplicial $B_n$-module. Since $A$ is noetherian,
we can apply \cite{Quillen} (Theorem~6.15) to conclude that for each
$q, n_0 \ge 0$, there exists an $N$ such that the natural map 
$D_q({B_{nn_0}}/A, B_{n_0}) \to D_q({B_{n_0}}/A, B_{n_0}) =
D_q({B_{n_0}}/A)$ is zero for all $n \ge N$. Since the map 
$D_q({B_{nn_0}}/A) \to D_q({B_{n_0}}/A)$ is the composite of the map 
\[D_q({B_{nn_0}}/A) = D_q({B_{nn_0}}/A, B_{nn_0}) \to
D_q({B_{nn_0}}/A, B_{n_0}) \to D_q({B_{n_0}}/A, B_{n_0}),\]   
we see that the map $D_q({B_{nn_0}}/A) \to D_q({B_{n_0}}/A)$ is zero
for all $n \ge N$. This in turn implies that the natural map
$D_q({B_{n}}/A) \to D_q({B_{n_0}}/A)$ is zero for all $n \gg n_0$
(in fact for all $n \ge Nn_0$). Now we apply Corollary~\ref{cor:exterior}
to the inverse system $B_{\bullet} = \{B_n\}$ and the free simplicial module
$M^{\bullet}_* = \{{\mathbb L}_{{B_n}/A}\}$ to conclude the proof of the
corollary.
\end{proof}
\section{Artin-Rees theorem for Andr{\'e}-Quillen Homology}
Let $k$ be a field and let $A$ be a $k$-algebra which is essentially of
finite type over $k$. Let $I$ be an ideal of $A$ and put 
$B_n = A/{I^{n+1}}$ for $n \ge 0$. Then $B_{\bullet} = \{B_n\}$ is an
inverse system of finite $A$-algebras. Moreover, ${\left \{
{\frac {D^{(r)}_q({B_{n}}/k)} {D^{(r)}_q({A}/k)}} \right \}}_{n \ge 0}$
is a $B_{\bullet}$-module. Similarly, ${\left \{{\rm Ker}
\left(D^{(r)}_q({A}/k, B_n) \to 
D^{(r)}_q({B_{n}}/k)\right)\right \}}_{n \ge 0}$ is also a 
$B_{\bullet}$-module. Our aim in this section is to prove an Artin-Rees
type theorem for these two modules. 
We begin with the following elementary result.   
\begin{lem}\label{lem:exact}
Let $A$ be any $k$-algebra and let 
\[ 0 \to M'_* \to M_* \to M''_* \to 0\]
be a short exact sequence of free simplicial $A$-modules. Then there
exists a convergent spectral sequence
\[ E^1_{p,q} = H_{q-p}\left({\wedge}^p M'_* {\otimes}_A 
{\wedge}^{r-p} M''_*\right) \Rightarrow  H_{q-p}\left({\wedge}^r M_*\right).\]
This spectral sequence is natural for morphisms of $k$-algebras and
morphisms of short exact sequences of free simplicial modules.
\end{lem}
\begin{proof} Exactness of simplicial modules means that it is exact at
each level and the exactness is compatible with the face and the degeneracy 
maps. For each $i \ge 0$, we can define a decreasing finite filtration on 
${\wedge}^r M_i$ by defining $F^j {\wedge}^r M_i$ to be the $A$-submodule
generated by the forms of the type
\[\left\{a_1 \wedge \cdots \wedge a_r | a_{i_1}, \cdots , a_{i_j} \in
M'_i \ {\rm for \ some} \ 1 \le i_1 \le \cdots \le i_j \le r \right\}.\]
Then we have
\[{\wedge}^r M_i = F^0 {\wedge}^r M_i \supseteq \cdots \supseteq
F^r {\wedge}^r M_i \supseteq F^{r+1} {\wedge}^r M_i = 0\] and it is easy to 
check that  for $0 \le j \le r$, the map
\[{\beta}^j_i : {\wedge}^j M'_i {\otimes}_A {\wedge}^{r-j} M_i \to
F^j {\wedge}^r M_i,\]
\[{\beta}^j_i \left((a_1 \wedge \cdots \wedge a_j) \otimes
(b_1 \wedge \cdots \wedge b_{r-j})\right) =
a_1 \wedge \cdots \wedge a_r \wedge b_1 \wedge \cdots \wedge b_{r-j}\]
descends to an isomorphism of quotients
\begin{equation}\label{eqn:descend} 
{\beta}^j_i :
{\wedge}^j M'_i {\otimes}_A {\wedge}^{r-j} M''_i \xrightarrow {\cong}
\frac {F^j {\wedge}^r M_i} {F^{j+1} {\wedge}^r M_i}.
\end{equation}
We also see from the above definition of the filtration and the maps
${\beta}^j_i$ that this filtration and the isomorphisms in 
~\ref{eqn:descend} are compatible with the morphisms of short exact
sequences. In particular, they are compatible with the face and the
degeneracy maps. Thus we get a decreasing filtration ${\left\{
F^j {\wedge}^r M_* \right\}}_{0 \le j \le r}$ of the simplicial module
${\wedge}^r M_*$  such that for each $0 \le j \le r$, there is a natural 
isomorphism
\begin{equation}\label{eqn:descend1} 
{\wedge}^j M'_* {\otimes}_A {\wedge}^{r-j} M''_* \xrightarrow {\cong}
\frac {F^j {\wedge}^r M_*} {F^{j+1} {\wedge}^r M_*}.
\end{equation}  
This filtration on the simplicial module ${\wedge}^r M_*$ gives 
(\cite{Weibel1}, 5.5) a convergent spectral sequence
\[E^1_{p, q} = H_{q-p}\left(\frac {F^p {\wedge}^r M_*} 
{F^{p+1} {\wedge}^r M_*}\right) \Rightarrow 
H_{q-p}\left({\wedge}^r M_*\right)\] 
with differential $E^1_{p, q} \to E^1_{p+1, q}$. The isomorphism of
~\ref{eqn:descend1} now completes the proof of the existence of
the spectral sequence. The functoriality with the morphisms of 
$k$-algebras and morphisms of exact sequences of simplicial modules
is clear from the definition of the filtration above, which is 
preserved under a morphism of exact sequences.
\end{proof}
\begin{cor}\label{cor:exact1} 
Let $A$ be an essentially of finite type algebra over a field $k$ and let 
$l \subset k$ be any subfield. Then there is a convergent spectral
sequence
\[E^1_{p, q} = {\Omega}^p_{k/l} {\otimes}_k D^{(r-p)}_{q-p}({A}/k)
\Rightarrow  D^{(r)}_{q-p}({A}/l).\]
\end{cor}
\begin{proof} We have (\cite{Quillen}, proof of Theorem~5.1) a short
exact sequence of free simplicial $A$-modules
\[ 0 \to {\mathbb L}_{k/l} {\otimes}_k A \to {\mathbb L}_{A/l} \to
{\mathbb L}_{A/k} \to 0.\]
Put ${\mathbb K}_{A/l} = {\mathbb L}_{k/l} {\otimes}_k A$. Then 
Lemma~\ref{lem:exact} gives us a convergent spectral sequence
\begin{equation}\label{eqn:E1}
E^1_{p,q} = H_{q-p}\left({\wedge}^p {\mathbb K}_{A/l} {\otimes}_A 
{\wedge}^{r-p} {\mathbb L}_{A/k} \right) \Rightarrow  H_{q-p}
\left({\wedge}^r  {\mathbb L}_{A/l} \right).
\end{equation}
To identify the $E^1$-terms, we see from the proof of 
Proposition~\ref{prop:tensor} that for each $p, q \ge 0$, there is a 
convergent spectral sequence
\['E^2_{i,j} = H_{i} \left(H_j\left({\wedge}^p {\mathbb K}_{A/l}\right)
{\otimes}_A {\wedge}^{r-p} {\mathbb L}_{A/k} \right) \Rightarrow 
H_{i+j}\left({\wedge}^p {\mathbb K}_{A/l} {\otimes}_A 
{\wedge}^{r-p} {\mathbb L}_{A/k} \right).\]
Since $A$ is $k$-flat, we have 
\[H_j\left({\wedge}^p {\mathbb K}_{A/l}\right) =
H_j\left({\wedge}^p {\mathbb L}_{k/l} {\otimes}_k A \right)    
= H_j\left({\wedge}^p {\mathbb L}_{k/l}\right) {\otimes}_k A
= D^{(p)}_{j}({k}/l) {\otimes}_k A.\]
Since this last group is a free $A$-module, we obtain
\['E^2_{i,j} =  D^{(p)}_{j}({k}/l) {\otimes}_k D^{(r-p)}_{i}({A}/k).\]
Now as $k$ is a direct limit of its subfields which are finitely
generated over $l$, and since the Andr{\'e}-Quillen homology commutes
with direct limits (\cite{Quillen}, 4.11), we see that $D^{(p)}_{j}({k}/l)$
is a direct limit of the Andr{\'e}-Quillen homology of the subfields
of $k$ which are finitely generated over $l$. In particular, we have
(\cite{Loday}, Theorem~3.5.6) 
\[D^{(p)}_{j}({k}/l) = \left\{ \begin{array}{ll}
{\Omega}^p_{k/l} & \mbox{if $p \ge 0, j = 0$} \\
0 & \mbox {otherwise} \hspace*{1cm}.
\end{array}
\right. \]
Thus we get
\['E^2_{i,j} = \left\{ \begin{array}{ll}
{{\Omega}^p_{k/l}{\otimes}_k D^{(r-p)}_{i}({A}/k)} &
\mbox{if $i \ge 0, j = 0$} \\
0 & \mbox{if $j > 0$} \hspace*{1cm}. 
\end{array}
\right. \]  
In particular, this spectral sequence degenerates at $'E^2$ and we get for
$p, i \ge 0$,
\[H_{i}\left({\wedge}^p {\mathbb K}_{A/l} {\otimes}_A 
{\wedge}^{r-p} {\mathbb L}_{A/k} \right) =
{\Omega}^p_{k/l} {\otimes}_k D^{(r-p)}_{i}({A}/k).\]
Putting this in our spectral sequence of ~\ref{eqn:E1}, we get the proof
of the corollary.
\end{proof}
Let $A$ be an essentially of finite type algebra over a field $k$ and
let $B_{\bullet} = \{B_n = A/{I^{n+1}}\}$ be the inverse system of
finite $A$-algebras as defined in the beginning of this section.  
\begin{lem}\label{lem:coeff}
For any given $r, q, n_0 \ge 0$, the natural map
\[\frac {D^{(r)}_{q}({A}/k, B_n)}{D^{(r)}_{q}({A}/k)} \to
\frac {D^{(r)}_{q}({A}/k, B_{n_0})}{D^{(r)}_{q}({A}/k)}\]
is zero for all $n \gg n_0$.
\end{lem}
\begin{proof}
For $r = 0$, both sides are zero, so we can assume $r \ge 1$.
We first observe that for any $n \ge 0$, one has
\[{\frac {D^{(r)}_{q}({A}/k, B_n)}{D^{(r)}_{q}({A}/k)}} 
\xrightarrow {\cong} 
{\frac {D^{(r)}_{q}({A}/k, B_n)}{D^{(r)}_{q}({A}/k) {\otimes}_A B_n}}.\]
By \cite{Quillen} (4.7), there is a convergent spectral sequence
\[{^n}E^2_{p, q} = {Tor}^A_p\left(D^{(r)}_{q}({A}/k), B_n\right)
\Rightarrow  D^{(r)}_{p+q}({A}/k, B_n).\]
This spectral sequence is compatible with the maps $B_n \surj B_{n-1}$
and gives a finite filtration of $D^{(r)}_{q}({A}/k, B_n)$
\[0 = F^n_{-1} \subseteq F^n_{0} \subseteq \cdots \subseteq F^n_{q-1} 
\subseteq F^n_{q} = D^{(r)}_{q}({A}/k, B_n)\]
such that ${^n}E^{\infty}_{j, q-j} = {F^n_{j}}/{F^n_{j-1}}$
for $0 \le j \le q$ and the edge map gives
${D^{(r)}_{q}({A}/k) {\otimes}_A B_n} \surj F^n_{0}$. Hence it suffices to
show that the natural map
\begin{equation}\label{eqn:coeff1}
\frac {D^{(r)}_{q}({A}/k, B_n)}{F^n_{0} D^{(r)}_{q}({A}/k, B_n)}  
\to \frac {D^{(r)}_{q}({A}/k, B_{n_0})}{F^{n_0}_{0} 
D^{(r)}_{q}({A}/k, B_{n_0})}
\end{equation}
is zero for all $n \gg n_0$.
Using the above filtration, an induction on $j$ and the doubling trick
of ~\ref{eqn:doubling}, this is reduced to showing that for $1 \le j \le q$,  
the map ${^n}E^2_{j, q-j} \to {^{n_0}}E^2_{j, q-j}$ is zero for all
$n \gg n_0$. But for $j \ge 1$, we have ${^n}E^2_{j, q-j} = 
{Tor}^A_j\left(D^{(r)}_{q-j}({A}/k), B_n\right)$. Furthermore, $A$ is 
a localization of a finite type $k$-algebra and so by \cite{Quillen}
(Proposition~4.12, Theorem~5.4(i)), $D^{(r)}_{q}({A}/k)$ is a finite
$A$-module for all $r, q \ge 0$. Hence by \cite{Andre}
(Proposition~10, Lemma~11), the map 
\[{Tor}^A_j\left(D^{(r)}_{q-j}({A}/k), B_n\right) \to
{Tor}^A_j\left(D^{(r)}_{q-j}({A}/k), B_{n_0}\right)\]
is zero for all $n \gg n_0$.
\end{proof}     
\begin{prop}\label{prop:coeff2}
Let $A$ and $B_{\bullet}$ be as above. Then for any given $r, q, n_0 \ge 0$,
the natural map
\[{\frac {D^{(r)}_{q}({B_n}/k)}{D^{(r)}_{q}({A}/k)}} \to
{\frac {D^{(r)}_{q}({B_{n_0}}/k)}{D^{(r)}_{q}({A}/k)}}\] is zero for
all $n \gg n_0$.
\end{prop}
\begin{proof} For $r= 0$, both sides are zero, so we assume $r \ge 1$.
Since the map $D^{(r)}_{q}({A}/k) \to 
D^{(r)}_{q}({B_n}/k)$ factors through the map
$D^{(r)}_{q}({A}/k) \to D^{(r)}_{q}({A}/k, B_n)$, one has for all $n$,
the natural exact sequence
\begin{equation}\label{eqn:coeff3}
{\frac {D^{(r)}_{q}({A}/k, B_n)}{D^{(r)}_{q}({A}/k)}} \to
{\frac {D^{(r)}_{q}({B_n}/k)}{D^{(r)}_{q}({A}/k)}} \to
{\frac {D^{(r)}_{q}({B_n}/k)}{D^{(r)}_{q}({A}/k, B_n)}} \to 0.
\end{equation}
Using the doubling trick of ~\ref{eqn:doubling} and Lemma~\ref{lem:coeff},
we only need to show that for any given $r \ge 1$ and $q, n_0 \ge 0$, the 
natural map
\begin{equation}\label{eqn:coeff4}
{\frac {D^{(r)}_{q}({B_n}/k)}{D^{(r)}_{q}({A}/k, B_n)}} \to
{\frac {D^{(r)}_{q}({B_{n_0}}/k)}{D^{(r)}_{q}({A}/k, B_{n_0})}}
\end{equation}
is zero for all $n \gg n_0$. 

For $n \ge 0$, we put 
${\mathbb K}_{{B_n}/k} = {\mathbb L}_{A/k} {\otimes}_A B_n$. Then we 
observe that for $r , q \ge 0$, $D^{(r)}_{q}({A}/k, B_n)$ is same
as $H_{q}\left(\left({\wedge}^r_A {\mathbb L}_{A/k}\right) {\otimes}_A
B_n \right) = H_{q}\left({\wedge}^r_{B_n} {\mathbb K}_{{B_n}/k}\right)$.
For all $n \ge 0$, we have an exact sequence (\cite{Quillen}, Theorem~5.1)
of free simplicial $B_n$-modules
\[ 0 \to {\mathbb L}_{A/k} {\otimes}_A B_n \to {\mathbb L}_{{B_n}/k} \to
{\mathbb L}_{{B_n}/A} \to 0.\]
Hence by Lemma~\ref{lem:exact}, there is a convergent spectral sequence
\[{^n}E^1_{p,q} = H_{q-p}\left({\wedge}^p {\mathbb K}_{{B_n}/k} 
{\otimes}_{B_n} {\wedge}^{r-p} {\mathbb L}_{{B_n}/A}\right) \Rightarrow  
H_{q-p}\left({\wedge}^r {\mathbb L}_{{B_n}/k}\right).\]
This spectral sequence is compatible with the maps $B_n \surj B_{n-1}$ and
gives a finite filtration of $H_{q}\left({\wedge}^r 
{\mathbb L}_{{B_n}/k}\right)$
\begin{equation}\label{eqn:coeff5}
H_{q}\left({\wedge}^r {\mathbb L}_{{B_n}/k}\right) = F^n_0 \supseteq
F^n_1 \supseteq \cdots \supseteq F^n_{r} \supseteq F^n_{r+1} = 0 
\end{equation}
with ${F^n_{j}}/{F^n_{j+1}} \cong {^n}E^{\infty}_{j,q+j}$ for $0 \le j \le r$
and a morphism of filtered modules 
$H_{q}\left({\wedge}^r {\mathbb L}_{{B_n}/k}\right)
\to H_{q}\left({\wedge}^r {\mathbb L}_{{B_{n-1}}/k}\right)$.
Furthermore, the edge map gives a surjection 
$H_{q}\left({\wedge}^r {\mathbb K}_{{B_n}/k}\right) \surj
F^n_r H_{q}\left({\wedge}^r {\mathbb L}_{{B_n}/k}\right)$.
In particular, we have 
${\frac {H_{q}\left({\wedge}^r {\mathbb L}_{{B_n}/k}\right)}
{H_{q}\left({\wedge}^r {\mathbb K}_{{B_n}/k}\right)}} \xrightarrow
{\cong} {\frac {H_{q}\left({\wedge}^r {\mathbb L}_{{B_n}/k}\right)}
{F^n_r H_{q}\left({\wedge}^r {\mathbb L}_{{B_n}/k}\right)}}$.
Hence using the filtration in ~\ref{eqn:coeff5}, an induction on $j$
and the doubling trick of ~\ref{eqn:doubling}, it suffices to show that
for $0 \le j \le r-1$, the natural map ${^n}E^1_{j, q+j} \to 
{^{n_0}}E^1_{j, q+j}$ is zero for all $n \gg n_0$.
But this follows from Proposition~\ref{prop:tensor} and 
Corollary~\ref{cor:AQ1}. This proves ~\ref{eqn:coeff4} and hence the
proposition.
\end{proof}
\begin{lem}\label{lem:kercoeff}
Let the $k$-algebras $A$ and $B_{\bullet}$ be as above. Then for any 
given $r, q, n_0 \ge 0$, the natural map
\[{\rm Ker}\left(D^{(r)}_{q}({A}/k, B_n) \to D^{(r)}_{q}({B_n}/k)\right)
\to {\rm Ker}\left(D^{(r)}_{q}({A}/k, B_{n_0}) \to 
D^{(r)}_{q}({B_{n_0}}/k)\right)\]
is zero for all $n \gg n_0$.
\end{lem}
\begin{proof} For $r = 0$, both sides are zero, so we assume $r \ge 1$.
We have seen in the proof of Proposition~\ref{prop:coeff2} that 
for all $n \ge 0$, $D^{(r)}_{q}({B_n}/k) = 
H_q\left({\wedge}^r {\mathbb L}_{{B_n}/k}\right)$ has a finite filtration
${\{F^n_j\}}_{0 \le j \le r}$ such that 
$H_{q}\left({\wedge}^r {\mathbb K}_{{B_n}/k}\right) \surj
F^n_r H_{q}\left({\wedge}^r {\mathbb L}_{{B_n}/k}\right)$. Thus we can 
replace $D^{(r)}_{q}({B_n}/k)$ by 
$F^n_r H_{q}\left({\wedge}^r {\mathbb L}_{{B_n}/k}\right) =
{^n}E^{\infty}_{r, q+r} = {^n}E^{r+1}_{r, q+r}$ in the statement of the lemma.

Now for $0 \le j \le r$, there is an exact sequence 
\[{^n}E^{j+1}_{r-j-1, q+r+j} \to {^n}E^{j+1}_{r, q+r} \to 
{^n}E^{j+2}_{r, q+r} \to 0,\]
which gives a finite increasing filtration of 
${\rm Ker}\left(D^{(r)}_{q}({A}/k, B_n) = {^n}E^{1}_{r, q+r}
\to {^n}E^{r+1}_{r, q+r}\right)$
\[0 = {\ov F}^n_{-1} \subseteq {\ov F}^n_{0} \subseteq \cdots \subseteq 
{\ov F}^n_{r-1} \subseteq {\ov F}^n_{r} = 
{\rm Ker}\left({^n}E^{1}_{r, q+r} \to {^n}E^{r+1}_{r, q+r}\right)\]
such that for all $0 \le j \le r$,
\[{^n}E^{j+1}_{r-j-1, q+r+j} \surj {{\ov F}^n_{j}}/{{\ov F}^n_{j-1}}.\] 
Thus by using an induction on $j$ and the doubling trick as before, it
suffices to show that for $0 \le j \le r$ and for given $n_0 \ge 0$,
the natural map
${^n}E^{1}_{r-j-1, q+r+j} \to {^{n_0}}E^{1}_{r-j-1, q+r+j}$ is zero
for all $n \gg n_0$. But this again follows from 
Proposition~\ref{prop:tensor} and Corollary~\ref{cor:AQ1}.
\end{proof}
\section{Andr{\'e}-Quillen Homology and Normalization I}
Let $A$ be an integral domain which is essentially of finite type algebra
over a field $k$. Assume that $A$ is singular, and let $f : A \to B$ be the 
normalization morphism of $A$.
We assume that $B$ is smooth over $k$. Our aim in this section is to
estimate the Andr{\'e}-Quillen homology of the conducting ideals for this
normalization. We begin with estimating the kernels of the maps
between differential forms. For any conducting ideal $I \subset A$ for the
normalization and for $r \ge 0$, let ${\Omega}^r_{(A, I)/l}$
(resp ${\Omega}^r_{(B, I)/l}$) denote the kernel of the map
${\Omega}^r_{A/l} \surj {\Omega}^r_{{(A/I)}/l}$ (resp
${\Omega}^r_{B/l} \surj {\Omega}^r_{{(B/I)}/l}$) for any subfield 
$l \subset k$.
\begin{lem}\label{lem:differential}
Let $A$ and $B$ be as above. Then for any given conducting ideal $I \subset
A$ for the normalization and for any $r \ge 0$, the natural map
\[{\Omega}^r_{(A, I^n)/k} \to {\Omega}^r_{(B, I^n)/k}\]
is injective for all sufficiently large $n$.
\end{lem}
\begin{proof} We see from ~\ref{eqn:D0} and Lemma~\ref{lem:kercoeff}
that for any given $n_0 \ge 0$ the natural map
\begin{equation}\label{eqn:D1}
{\rm Ker} \left({\Omega}^r_{A/k} {\otimes}_A A/{I^n} \xrightarrow {u^A_n}
{\Omega}^r_{{(A/{I^n})}/k}\right) \to
{\rm Ker} \left({\Omega}^r_{A/k} {\otimes}_A A/{I^{n_0}} \xrightarrow 
{u^A_{n_0}} {\Omega}^r_{{(A/{I^{n_0}})}/k}\right)
\end{equation}
is zero for all $n \gg n_0$. In the same way, the map
${\rm Ker}(u^B_n) \to {\rm Ker}(u^B_{n_0})$ is zero for all $n \gg n_0$. 
On the other hand, we have a commutative diagram of exact sequences 
for any $n \ge 0$.
\[
\xymatrix{
0 \ar[r] & I^n {\Omega}^r_{A/k} \ar[r] \ar[d] & {\Omega}^r_{(A, I^n)/k}
\ar[d] \ar[r] & {\rm Ker}(u^A_n) \ar[r] \ar[d] & 0 \\
0 \ar[r] & I^n {\Omega}^r_{B/k} \ar[r] & {\Omega}^r_{(B, I^n)/k}
\ar[r] & {\rm Ker}(u^B_n) \ar[r] & 0 }
\]
We claim that the map $I^n {\Omega}^r_{A/k} \to I^n {\Omega}^r_{B/k}$ is
injective for all $n \gg 0$.
To see this, note that ${\Omega}^r_{A/k}$ is a finite $A$-module
and we can apply the Artin-Rees theorem to find a $c > 0$ such that
all $n > c$, one has 
\[\left(I^n {\Omega}^r_{A/k} \cap {\Omega}^r_{(A, B)/k}\right)
\subseteq I^{n-c}\left(I^c {\Omega}^r_{A/k} \cap {\Omega}^r_{(A, B)/k}\right)
\subseteq I^{n-c} {\Omega}^r_{(A, B)/k},\]
where ${\Omega}^r_{(A, B)/k} = {\rm Ker}\left({\Omega}^r_{A/k} \to
{\Omega}^r_{B/k}\right)$.
On the other hand, the finite $A$-module ${\Omega}^r_{(A, B)/k}$ is 
supported on the support of $I$ and hence $I^{n-c} {\Omega}^r_{(A, B)/k}
= 0$ for all $n \gg 0$. This proves the claim.   
Using the claim in the above diagram, we see that there exists an $n_0$
such that ${\rm Ker}\left({\Omega}^r_{(A, I^n)/k} \to 
{\Omega}^r_{(B, I^n)/k}\right) \inj {\rm Ker}(u^A_n)$ for all $n \ge n_0$.
Now we use ~\ref{eqn:D1} to conclude that the map
\[{\rm Ker}\left({\Omega}^r_{(A, I^n)/k} \to 
{\Omega}^r_{(B, I^n)/k}\right) \to 
{\rm Ker}\left({\Omega}^r_{(A, I^{n_0})/k} \to 
{\Omega}^r_{(B, I^{n_0})/k}\right)\]
is zero for all $n \gg n_0$. However, this map is clearly injective.
Hence we must have ${\rm Ker}\left({\Omega}^r_{(A, I^n)/k} \to 
{\Omega}^r_{(B, I^n)/k}\right) = 0$ for all $n \gg 0$.
\end{proof}
\begin{lem}\label{lem:AQN}
Let $f: A \to B$ be as above. Then for any conducting ideal $I$ for the
normalization and for any $r, q \ge 1$, the natural map
\[D^{(r)}_{q}({A}/k) \xrightarrow {u^A_n} D^{(r)}_{q}({A/{I^n}}/k)\]
is injective for all $n \gg 0$.
\end{lem}
\begin{proof} We first observe that
$D^{(r)}_{q}({B}/k) = 0$ for $q \ge 1$ as $B$ is
smooth (\cite{Loday}, Theorem~3.5.6). Since the homology groups
$D^{(r)}_{q}({A}/k)$ are finite $A$-modules and since they
commute with the localization (\cite{Quillen}, Theorem~5.4 (i)), we
see as before that $I^n D^{(r)}_{q}({A}/k) = 0$ for all $n \gg 0$
whenever $q \ge 1$. In particular, for any $r, q \ge 1$
there exits $N \gg 0$ such that for all $n \ge N$, 
\begin{equation}\label{eqn:AQN1}
D^{(r)}_{q}({A}/k) \xrightarrow {\cong} D^{(r)}_{q}({A}/k) {\otimes}_A
A/{I^n}.
\end{equation}

For $r, q \ge 1$ and $n \ge 0$, we have a convergent spectral sequence
\[{^n}E^2_{p, q} = {Tor}^A_p\left(D^{(r)}_{q}({A}/k), A/{I^n}\right)
\Rightarrow  D^{(r)}_{p+q}({A}/k, A/{I^n})\]
with differential ${^n}E^2_{p, q} \to {^n}E^2_{p-2, q+1}$. This gives
a filtration 
\[0 = F^n_{-1} \subseteq F^n_{0} \subseteq \cdots \subseteq 
F^n_{q-1} \subseteq F^n_{q} =  D^{(r)}_{q}({A}/k, A/{I^n})\]
and map of filtered modules $D^{(r)}_{q}({A}/k, A/{I^n}) \to
D^{(r)}_{q}({A}/k, A/{I^{n-1}})$. The edge map further gives
a surjection $D^{(r)}_{q}({A}/k) {\otimes}_A A/{I^n} 
= {^n}E^2_{0, q} \surj F^n_{0} = {^n}E^{\infty}_{0, q} = 
{^n}E^{q+2}_{0, q}$.  

We now show that for $r, q \ge 1$ and $n_0 \ge 0$, the natural map
\begin{equation}\label{eqn:AQN2}
{\rm Ker}\left(D^{(r)}_{q}({A}/k) {\otimes}_A A/{I^n} \to
 F^n_{0}  D^{(r)}_{q}({A}/k, A/{I^n})\right) \hspace*{4cm}
\end{equation}
\[\hspace*{5cm} \to
{\rm Ker}\left(D^{(r)}_{q}({A}/k) {\otimes}_A A/{I^{n_0}} \to
 F^{n_0}_{0}  D^{(r)}_{q}({A}/k, A/{I^{n_0}})\right)\]
is zero for all $n \gg n_0$.

For $0 \le j \le q$, there is an exact sequence
\[{^n}E^{j+2}_{j+2, q-j-1} \to {^n}E^{j+2}_{0, q} \to {^n}E^{j+3}_{0, q}
\to 0.\]
Thus by letting ${\Gamma}^n_j = {\rm Ker} \left({^n}E^2_{0, q} \surj
{^n}E^{j+2}_{0, q}\right)$ for $0 \le j \le q$, we get a filtration
$0 = {\Gamma}^n_0 \subseteq {\Gamma}^n_1 \subseteq \cdots \subseteq
{\Gamma}^n_q$ of ${\Gamma}^n_q$ such that
\[{^n}E^{j+2}_{j+2, q-j-1} \surj {{\Gamma}^n_{j+1}}/{{\Gamma}^n_j}.\]
Thus to prove ~\ref{eqn:AQN2}, it suffices to show by an induction
on $j$ and the doubling trick that for given $r, q, n_0 \ge 1$ and
for $0 \le j \le q$, the natural map ${^n}E^{2}_{j+2, q-j-1}
\to {^{n_0}}E^{2}_{j+2, q-j-1}$ is zero for $n \gg n_0$.
But for $n, j \ge 0$ we have ${^n}E^{2}_{j+2, q-j-1}= $ \\
${Tor}^A_{j+2}\left(D^{(r)}_{q-j-1}({A}/k), A/{I^n}\right)$ and the
map \[{Tor}^A_{j+2}\left(D^{(r)}_{q-j-1}({A}/k), A/{I^n}\right)
\to {Tor}^A_{j+2}\left(D^{(r)}_{q-j-1}({A}/k), A/{I^{n_0}}\right)\]
is zero by for $n \gg n_0$ by \cite{Andre}(Proposition~10, Lemma~11) as 
$D^{(r)}_{q}({A}/k)$ are all finite $A$-modules. This proves 
~\ref{eqn:AQN2}.

If $D^{(r)}_{q}({A}/k) \xrightarrow {{\theta}^A_n} 
D^{(r)}_{q}({A}/k, A/{I^n})$ and $D^{(r)}_{q}({A}/k, A/{I^n})
\xrightarrow {v^A_n} D^{(r)}_{q}({A/{I^n}}/k)$ denote the natural maps,
then for any $r, q, n \ge 1$, we get a natural exact sequence
\[0 \to {\rm Ker}\left({{\theta}^A_n}\right) \to
{\rm Ker}\left(u^A_n\right) \to {\rm Ker}\left(v^A_n\right).\]
Now for any $n_0 \ge N$, we can use ~\ref{eqn:AQN1} and  ~\ref{eqn:AQN2}
to conclude that the map ${\rm Ker}\left({{\theta}^A_n}\right) 
\to {\rm Ker}\left({{\theta}^A_{n_0}}\right)$ is zero for
all $n \gg n_0$. The map ${\rm Ker}\left(v^A_n\right) \to
{\rm Ker}\left(v^A_{n_0}\right)$ is zero for all $n \gg n_0$ by
Lemma~\ref{lem:kercoeff}. The doubling trick again shows that for
any $r, q \ge 1$ and $n_0 \ge N$, the natural map
${\rm Ker}\left(u^A_n\right) \to {\rm Ker}\left(u^A_{n_0}\right)$
is zero for all $n \gg n_0$. But this last map is clearly injective.
Hence we must have ${\rm Ker}\left(u^A_n\right) = 0$ for all $n \gg 0$.
\end{proof}  
\section{Andr{\'e}-Quillen Homology and Normalization II}
Most of the proofs in the previous section relied on the fact that
the algebra $A$ is essentially of {\sl finite type} over the field $k$.
Our aim in this section is to generalize the results of the previous
section to the case when the base ring for the Andr{\'e}-Quillen
homology of $k$-algebras is any subfield of $k$. Our eventual application
will need these results when the base ring is the field of rational numbers.
So let $A$ be an integral domain which is essentially of finite type
algebra over a field $k$. Let $f : A \to B$ be the normalization of
$A$ such that $B$ is smooth. Let $l \subset k$ be any subfield.

The basic extra ingredient to deal with the general case will be our 
spectral sequence of Corollary~\ref{cor:exact1}:
\begin{equation}\label{eqn:generalfield}
E^1_{p, q}(A) = {\Omega}^p_{k/l} {\otimes}_k D^{(r-p)}_{q-p}({A}/k)
\Rightarrow  D^{(r)}_{q-p}({A}/l).
\end{equation}
As shown in Lemma~\ref{lem:exact}, this spectral sequence is clearly 
compatible with the morphisms of 
$k$-algebras. We denote the corresponding spectral sequence
for $A/{I^n}$ by ${^n}E^i_{p, q}$ as before and that for $B$ by
$E^i_{p, q}(B)$. Put $E^i_{p, q}(A,B)= {\rm Ker}\left(E^i_{p, q}(A)
\to E^i_{p, q}(B)\right)$.
\begin{lem}\label{lem:spectral}
For any $r, i \ge 1$ and $p, q, n_0 \ge 0$, the natural map
\begin{equation}\label{eqn:spectral1}
{\frac {{^n}E^i_{p, q}} {E^i_{p, q}(A)}} \to
{\frac {{^{n_0}}E^i_{p, q}} {E^i_{p, q}(A)}}
\end{equation}
is zero for all $n \gg n_0$. Furthermore, the natural map
\begin{equation}\label{eqn:spectral2}
E^i_{p, q}(A,B) \to {^n}E^i_{p, q}
\end{equation}
is injective for all $n \gg n_0$.
\end{lem}
\begin{proof} We prove both statements by induction on $i \ge 1$.
For $i =1$, ~\ref{eqn:spectral1} follows directly from 
Proposition~\ref{prop:coeff2} and ~\ref{eqn:spectral2} follows directly
from Lemmas ~\ref{lem:differential} and ~\ref{lem:AQN}.
So assume that  ~\ref{eqn:spectral1} and ~\ref{eqn:spectral2} hold for
all $1 \le j \le i$. We first show that ~\ref{eqn:spectral2} holds for
$i+1$. Consider the following commutative diagram of exact sequences.
\[
\xymatrix{
E^i_{p-i, q+i-1}(A) \ar[r] \ar[d] & 
{\rm Ker}\left({\del}^i_{p,q}(A)\right) \ar[d] 
\ar[r] & E^{i+1}_{p, q}(A) \ar[d] \ar[r] & 0 \\
E^i_{p-i, q+i-1}(B) \ar[r] & {\rm Ker}\left({\del}^i_{p,q}(B)\right)  
\ar[r] & E^{i+1}_{p, q}(B) \ar[r] & 0,}  
\]
where $E^{i}_{p, q} \xrightarrow {{\del}^i_{p,q}}
E^{i}_{p+i, q-i+1}$ is the differential of the spectral sequence.
If $p \neq q$, then $E^1_{p, q}(B) = 
{\Omega}^p_{k/l} {\otimes}_k D^{(r-p)}_{q-p}({B}/k) = 0$ (as $B$ is smooth
over $k$) and so is $E^j_{p, q}(B)$ for $j \ge 1$. This gives exact
sequence
\begin{equation}\label{eqn:spectral3}
E^i_{p-i, q+i-1}(A) \to {\rm Ker}\left({\rm Ker}\left({\del}^i_{p,q}(A)
\right) \to {\rm Ker}\left({\del}^i_{p,q}(B)\right)\right) \to
E^{i+1}_{p, q}(A, B) \to 0.
\end{equation}
If $p = q$, then $q+i-1-p+i = 2i-1 \ge 1$ and we get 
$E^i_{p-i, q+i-1}(B) = 0$. This again gives the exact sequence as above.
Thus ~\ref{eqn:spectral3} holds for all $p, q \ge 0$.   
Now we consider the following commutative diagram with exact rows.
\[
\xymatrix@C.5pc{
E^i_{p-i, q+i-1}(A) \ar[r] \ar[d] & 
{\rm Ker}\left({\rm Ker}\left({\del}^i_{p,q}(A)
\right) \to {\rm Ker}\left({\del}^i_{p,q}(B)\right)\right) \ar[r] \ar[d] & 
E^{i+1}_{p, q}(A, B)  \ar[r] \ar[d] & 0 \\
{^n}E^i_{p-i, q+i-1} \ar[r] & {\rm Ker}\left({^n}{\del}^i_{p,q}\right)
\ar[r] & {^n}E^{i+1}_{p, q} \ar[r] & 0}
\]
The middle vertical arrow is injective for $n \gg n_0$ by induction.
Let $N$ be the smallest integer such that this arrow is injective
for $n \ge N$. Then a diagram chase shows that 
${\rm Ker}\left(E^{i+1}_{p, q}(A, B) \to {^n}E^{i+1}_{p, q}\right)$
is an $A$-submodule of a quotient of 
${\frac {{^n}E^i_{p-i, q+i-1}}{E^i_{p-i, q+i-1}(A)}}$ for all $n \ge N$.
Thus to prove ~\ref{eqn:spectral2} for $i+1$, it suffices to show that
the natural map 
\[{\frac {{^n}E^i_{p-i, q+i-1}}{E^i_{p-i, q+i-1}(A)}} \to
{\frac {{^N}E^i_{p-i, q+i-1}}{E^i_{p-i, q+i-1}(A)}}\]
is zero for all $n \gg N$. But this is true as ~\ref{eqn:spectral1}
holds for $i$ by induction.   

Now we show that ~\ref{eqn:spectral1} holds for $i+1$.
We have a commutative diagram
\[
\xymatrix@C4.0pc{
E^{i}_{p, q}(A) \ar[r]^{{\del}^i_{p,q}(A)} 
\ar[d] & E^{i}_{p+i, q-i+1}(A) \ar[d] \\
E^{i}_{p, q}(B) \ar[r]^{{\del}^i_{p,q}(B)} & E^{i}_{p+i, q-i+1}(B).}
\]
If $p \neq q$, then $E^{1}_{p, q}(B) = 0 = E^{i}_{p, q}(B)$  and if
$p = q$, then $q-i+1-p-i = 1-2i < 0$ as $i \ge 1$ and hence
$E^{i}_{p+i, q-i+1}(B) = 0$. Now the above diagram shows that for
$p, q \ge 0$, one has a factorization
\begin{equation}\label{eqn:spectral4}
E^{i}_{p, q}(A) \xrightarrow {{\del}^i_{p,q}(A)} 
E^{i}_{p+i, q-i+1}(A,B) \inj E^{i}_{p+i, q-i+1}(A).
\end{equation}
Next we consider another commutative diagram for $n \ge 0$.
\[
\xymatrix{
0 \ar[r] & {\rm Ker}\left({\del}^i_{p,q}(A)\right) \ar[r] \ar[d] &
E^{i}_{p, q}(A) \ar[r] \ar[d] & {\rm Image}\left({\del}^i_{p,q}(A)\right)
\ar[r] \ar[d] & 0 \\
0 \ar[r] & {\rm Ker}\left({^n}{\del}^i_{p,q}\right) \ar[r] &
{^n}E^{i}_{p, q} \ar[r] & {\rm Image}\left({^n}{\del}^i_{p,q}\right)
\ar[r] & 0}
\]
Since ~\ref{eqn:spectral2} holds for $i$, we see from
~\ref{eqn:spectral4} that the right vertical
arrow is injective for all $n \gg n_0$. In particular, we get
an inclusion 
\[{\frac {{\rm Ker}\left({^n}{\del}^i_{p,q}\right)}
{{\rm Ker}\left({\del}^i_{p,q}(A)\right)}} \inj
{\frac {{^n}E^{i}_{p, q}}{E^{i}_{p, q}(A)}}\]
for all $n \ge N \gg n_0$. We now apply induction on $i$ in 
~\ref{eqn:spectral1} and this inclusion to conclude that the map
\[{\frac {{\rm Ker}\left({^n}{\del}^i_{p,q}\right)}
{{\rm Ker}\left({\del}^i_{p,q}(A)\right)}} \to
{\frac {{\rm Ker}\left({^N}{\del}^i_{p,q}\right)}
{{\rm Ker}\left({\del}^i_{p,q}(A)\right)}}\]
is zero for all $n \gg N$. Since 
${\frac {{^n}E^{i+1}_{p, q}}{E^{i+1}_{p, q}(A)}}$ is a quotient of 
the module on the left for $n \ge 0$, we get that the map
\[{\frac {{^n}E^{i+1}_{p, q}}{E^{i+1}_{p, q}(A)}} \to
{\frac {{^N}E^{i+1}_{p, q}}{E^{i+1}_{p, q}(A)}}\]
is zero for all $n \gg N$. Since $N \ge n_0$, we see that 
~\ref{eqn:spectral1} holds for $i \ge 1$. This proves the lemma.
\end{proof}
The following is the generalization of Proposition~\ref{prop:coeff2} 
and Lemmas ~\ref{lem:differential} and ~\ref{lem:AQN} to the case
when the base ring of the Andr{\'e}-Quillen homology is any subfield
of the given field $k$.
\begin{prop}\label{prop:generalcoeff}
Let $f : A \to B$ be the smooth normalization of an essentially
of finite type $k$-algebra $A$. Let $l \subset k$ be a subfield.
Let $I \subset A$ be any given conducting ideal for the normalization.
Then \\
$(i)$ For any $r, q, n_0 \ge 0$, the natural map
\[{\frac {D^{(r)}_{q}({A/{I^n}}/l)}{D^{(r)}_{q}({A}/l)}}
\to {\frac {D^{(r)}_{q}({A/{I^{n_0}}}/l)}{D^{(r)}_{q}({A}/l)}}\]
is zero for all $n \gg n_0$. \\
$(ii)$ For any $r, q, n_0 \ge 0$, the natural map
\[{\frac {D^{(r)}_{q}({B/{I^n}}/l)}{D^{(r)}_{q}({B}/l)}}
\to {\frac {D^{(r)}_{q}({B/{I^{n_0}}}/l)}{D^{(r)}_{q}({B}/l)}}\]
is zero for all $n \gg n_0$. \\ 
$(iii)$ For any $r, q \ge 1$, the natural maps
\[D^{(r)}_{q}({A}/l) \to D^{(r)}_{q}({A/{I^n}}/l)\]
\[D^{(r)}_{0}({(A,B)}/l) \to D^{(r)}_{q}({A/{I^n}}/l)\]
are injective for all $n \gg n_0$.
\end{prop}
\begin{proof} For $r =0$, the part $(i)$ is obvious as the groups on the
both sides are zero. So we assume $r \ge 1$.
The spectral sequence ~\ref{eqn:generalfield}
({\sl cf.} Corollary~\ref{cor:exact1}) gives for any $r \ge 1$ and 
$q \ge 0$, a finite filtration of $D^{(r)}_{q}({A}/l)$
\[D^{(r)}_{q}({A}/l) = F^0(A) \supseteq F^1(A) \supseteq \cdots \supseteq 
F^r(A) \supseteq F^{r+1}(A) = 0\]
with ${F^j(A)}/{F^{j+1}(A)} \cong E^{\infty}_{j, q+j}(A) = 
E^{r+1}_{j, q+j}(A)$ for $0 \le j \le r$.
One has similar filtrations for $D^{(r)}_{q}({B}/l)$ and
$D^{(r)}_{q}({A/{I^n}}/l)$ together with morphisms of filtered $A$-modules.
Thus for $0 \le j \le r$ and $n \ge 0$, we have the exact sequence
\[{\frac {F^{j+1}(A/{I^n})}{F^{j+1}(A)}} \to 
{\frac {F^{j}(A/{I^n})}{F^{j}(A)}} \to
{\frac {{^n}E^{r+1}_{j, q+j}}{E^{r+1}_{j, q+j}(A)}} \to 0.\]
Now by comparing this exact sequence for $n_0$ and $n \ge n_0$, using
the doubling trick in ~\ref{eqn:doubling} as before and using the 
descending induction on $j$, we see that it suffices to show that 
for $r \ge 1$, $q, n_0 \ge 0$ and $0 \le j \le r$, the natural map
\[{\frac {{^n}E^{r+1}_{j, q+j}}{E^{r+1}_{j, q+j}(A)}} \to
{\frac {{^{n_0}}E^{r+1}_{j, q+j}}{E^{r+1}_{j, q+j}(A)}}\]
is zero for all $n \gg n_0$ to prove part $(i)$ of the proposition.
But this follows from Lemma~\ref{lem:spectral}.

To prove $(ii)$, we first observe from the smoothness of $B$ and the
above spectral sequence that
\[D^{(r)}_{q}({B}/l) = \left\{ \begin{array}{ll}
{\Omega}^r_{B/l} & \mbox{if $q= 0$} \\
0 & \mbox {otherwise} \ \ .
\end{array}
\right. \]
Furthermore, as $D^{(r)}_{0}({B/{I^n}}/l) = {\Omega}^r_{{(B/{I^n})}/l}$
({\sl cf.} ~\ref{eqn:D0}), we immediately get $(ii)$ for $q = 0$.
We also get from this that for $q \ge 1$ and $n \ge 0$,
\begin{equation}\label{eqn:SS1}
{\frac {D^{(r)}_{q}({B/{I^n}}/l)}{D^{(r)}_{q}({B}/l)}}
\xrightarrow {\cong} D^{(r)}_{q}({B/{I^n}}/l).
\end{equation}
Thus we need to show that the map $D^{(r)}_{q}({B/{I^n}}/l) \to
D^{(r)}_{q}({B/{I^{n_0}}}/l)$ is zero for all $n \gg n_0$.
For $r, q \ge 1$, the natural map
\begin{equation}\label{eqn:SS} 
D^{(r)}_{q}({B/{I^n}}/k) \to
D^{(r)}_{q}({B/{I^{n_0}}}/k)
\end{equation} 
is zero for all $n \gg n_0$ by ~\ref{eqn:SS1} and
Proposition~\ref{prop:coeff2}. 
Now the proof of $(ii)$ follows by using ~\ref{eqn:SS} and the
spectral sequence ~\ref{eqn:generalfield} and then by following 
exactly the same argument as in the proof of $(i)$. Here the
analogue of Lemma~\ref{lem:spectral} follows immediately
from ~\ref{eqn:generalfield} and ~\ref{eqn:SS}.

For proving part $(iii)$, we consider the following
diagram of exact sequences for $0 \le j \le r$.
\[
\xymatrix{
0 \ar[r] & F^{j+1}(A,B) \ar[d] \ar[r] & F^{j}(A,B) \ar[r] \ar[d] &
E^{r+1}_{j, q+j}(A,B) \ar[d] & \\
0 \ar[r] & {^n}F^{j+1} \ar[r] & {^n}F^{j} \ar[r] &
{^n}E^{r+1}_{j, q+j} \ar[r] & 0, }
\]   
where $F^{j}(A,B) = {\rm Ker}\left(F^{j}(A) \to F^{j}(B)\right)$ and
${^n}F^{j}$ is the filtration of $D^{(r)}_{q}({A/{I^n}}/l)$.
Again by descending induction on $j$, it suffices to show that for 
$r \ge 1$, $q, n_0 \ge 0$ and for $0 \le j \le r$, the map
\[E^{r+1}_{j, q+j}(A,B) \to {^n}E^{r+1}_{j, q+j}\]
is injective. But this follows again from Lemma~\ref{lem:spectral}.
\end{proof}
\section{Andr{\'e}-Quillen to Hochschild Homology}
In this section, we derive some consequences of our results of the
previous section for the Hochschild and cyclic homology.
Let $k$ be a field and let $A$ be an essentially of finite
type $k$-algebra. We recall that for any $A$-algebra $B$, the
Hochschild homology $HH^A_*(B)$ are the homology group of the
pre-simplicial $B$-module $C^A_*(B)$ given by
\[C^A_n(B) = B {\otimes}_A \cdots {\otimes}_A B = B^{{\otimes}{n+1}} \
{\rm and \ the \ face\ maps}\]
$d_i : C^A_n(B) \to C^A_{n-1}(B)$ for $0 \le i \le n$ being given by
\[d_i(a_0, \cdots , a_n) = (a_0, \cdots , a_i a_{i+1}, \cdots , a_n)
\ {\rm for} \ 0 \le i \le n-1 \ {\rm and}\] 
\[d_n(a_0, \cdots , a_n) = (a_n a_0, a_1, \cdots, a_{n-1}).\]
The associated chain complex of $B$-modules is called the 
{\sl Hochschild} complex of $B$ over $A$.

To define the cyclic homology of $B$, one uses the natural action of
finite cyclic groups on the Hochschild complex to construct the
cyclic bicomplex $CC^A(B)$ and the cyclic homology $HC^A_*(B)$ of $B$ 
over $A$ are defined as the homology of the associated total complex.
We refer the reader to \cite{Loday} for the details about the
definitions of Hochschild and cyclic homology and their properties
which are relevant to us in this paper. For any map $B \to B'$ of
$A$-algebras, the relative Hochschild homology $HH^A_*(B, B')$ is
defined as the homology of the complex ${\rm Cone}
\left(HH^A_*(B) \to HH^A_*(B')\right)[-1]$. For an ideal $I \subset B$,
the relative Hochschild homology $HH^A_*(B, I)$ is the relative
homology of the map $B \to B/I$. For a map $B \to B'$ and an ideal
$I \subset B$ such that $IB' = I$, the double relative Hochshild
homology are defined as the homology of the complex \\
${\rm Cone}\left(HH^A_*(B, I) \to HH^A_*(B', I)\right)[-1]$.
The relative and double relative cyclic homology are defined in the
analogous way by taking the cones over the total cyclic complexes.
If ${CC^A(B)}^2$ denotes the cyclic bicomplex of $B$ consisting of 
only the first two columns of $CC^A(B)$, then there is a natural
short exact sequence ({\sl loc. cit.}, Theorem~2.2.1) 
\[0 \to {CC^A(B)}^2 \to CC^A(B) \to CC^A(B)[2, 0] \to 0,\]
which gives the Connes' periodicity long exact sequence 
(also called SBI-sequence)
\[ \cdots \to HH^A_n(B) \xrightarrow {I} HC^A_n(B) \xrightarrow {S}
HC^A_{n-2}(B) \xrightarrow {B}  HH^A_{n-1}(B) \xrightarrow {I} 
\cdots .\]
Taking the cones over the appropriate short exact sequences as above,
one gets similar SBI-sequence for the relative and double relative
Hochschild and cyclic homology.   

The most crucial fact which will be useful to us in this paper
is the canonical decomposition of the Hochschild homology in terms of
the Andr{\'e}-Quillen homology in characteristic zero.
We state it here for the sake of 
reader's convenience and refer to {\sl loc. cit.} (Theorem~3.5.8)
for the proof.
\begin{thm}\label{thm:AQH}
For any flat $A$-algebra $B$, there is a canonical decomposition
\[HH^A_n(B) \cong \underset {r+q = n} {\bigoplus} D^{(r)}_q(B/A).\]
\end{thm}
It is also known that this decomposition is compatible with the
Hodge decomposition of the Hochschild homology.

An immediate consequence of this canonical decomposition is the
following Hochschild homology analogue of 
Proposition~\ref{prop:generalcoeff}. As in the previous section,
let $A$ be an integral domain which is an essentially of finite
type algebra over a field $k$. Let $f : A \to B$ be the smooth
normalization of $B$. For any subfield $l \subset k$ and for
$i \in \Z$, we denote the kernel of the
map $HH^l_i(A) \to HH^l_i(B)$ by ${\ov {HH}}^l_i(A, B)$.
\begin{cor}\label{cor:AQH1}
Let $I \subset A$ be any given conducting ideal for the normalization
$f: A \to B$ as above. Then for any subfield $l \subset k$ and
any $i \ge 0$, \\
$(i)$ The natural map
\[{\frac {HH^l_i(A/{I^n})}{HH^l_i(A)}}
\to {\frac {HH^l_i(A/{I^{n_0}})}{HH^l_i(A)}}\]
is zero for all $n \gg n_0$. \\
$(ii)$ The natural map
\[{\frac {HH^l_i(B/{I^n})}{HH^l_i(B)}}
\to {\frac {HH^l_i(B/{I^{n_0}})}{HH^l_i(B)}}\]
is zero for all $n \gg n_0$. \\
$(iii)$ The natural map
\[{\ov {HH}}^l_i(A, B) \to HH^l_i(A/{I^n})\]
is injective for all $n \gg n_0$. 
\end{cor}
\begin{proof} This follows immediately from 
Proposition~\ref{prop:generalcoeff} and the canonical decomposition of the 
Hochschild homology in Theorem~\ref{thm:AQH}. 
\end{proof}
\begin{cor}\label{cor:AQH2}
Let the notations be as in Corollary~~\ref{cor:AQH1}. Then for any given
conducting ideal $I$ for the normalization and for any $i, n_0 \ge 0$,
the natural maps
\begin{equation}\label{eqn:AQH20}
{\frac {HH^l_i(B, I^n)}{HH^l_i(A, I^n)}} \to
{\frac {HH^l_i(B, I^{n_0})}{HH^l_i(A, I^{n_0})}}
\end{equation}
\begin{equation}\label{eqn:AQH21}
{\rm Ker}\left(HH^l_i(A, I^{n}) \to HH^l_i(B, I^{n})\right)
\to {\rm Ker}\left(HH^l_i(A, I^{n_0}) \to HH^l_i(B, I^{n_0})\right)
\end{equation}
are zero for all $n \gg n_0$.
\end{cor}
\begin{proof}  We consider the following commutative diagram of short 
exact sequences coming from the long exact relative Hochschild homology 
sequence.
\begin{equation}\label{eqn:AQH22}
\xymatrix@C.5pc{
0 \ar[r] & {\frac {HH^l_{i+1}(A/{I^n})}{HH^l_{i+1}(A)}} \ar[r] 
\ar[d] & HH^l_i(A, I^{n}) \ar[r] \ar[d] &
{\rm Ker}\left(HH^l_{i}(A) \to HH^l_{i}(A/{I^n})\right) \ar[r] \ar[d] & 0 \\
0 \ar[r] & {\frac {HH^l_{i+1}(B/{I^n})}{HH^l_{i+1}(B)}} \ar[r] 
& HH^l_i(B, I^{n}) \ar[r] &
{\rm Ker}\left(HH^l_{i}(B) \to HH^l_{i}(B/{I^n})\right) \ar[r] & 0}
\end{equation}
Using the naturality of the Hodge decomposition of the Hochschild
homology ({\sl loc. cit.}, Theorem~4.5.10) and smoothness of $B$, we have
the isomorphism 
\[{\frac {{\rm Ker}\left(HH^l_{i}(B) \to HH^l_{i}(B/{I^n})\right)}
{{\rm Ker}\left(HH^l_{i}(A) \to HH^l_{i}(A/{I^n})\right)}}
\xrightarrow {\cong} {\frac {{\Omega}^i_{(B, I^n)/l}}
{{\Omega}^i_{(A, I^n)/l}}}.\]
Using this identification, the above diagram gives us for $n \ge 0$, an
exact sequence of quotients
\[{\frac {HH^l_{i+1}(B/{I^n})}{HH^l_{i+1}(B)}} \to {\frac {HH^l_i(B, I^{n})}
{HH^l_i(A, I^{n})}} \to {\frac {{\Omega}^i_{(B, I^n)/l}}
{{\Omega}^i_{(A, I^n)/l}}} \to 0.\]
By \cite{Kr2} (Lemma~4.1), one has that for any given $i, n_0 \ge 0$,
the natural map ${\frac {{\Omega}^i_{(B, I^n)/l}}{{\Omega}^i_{(A, I^n)/l}}}
\to {\frac {{\Omega}^i_{(B, I^{n_0})/l}}{{\Omega}^i_{(A, I^{n_0})/l}}}$
is zero for all $n \gg n_0$. Corollary~\ref{cor:AQH1} implies that
for $i, n_0 \ge 0$, the natural map
${\frac {HH^l_{i}(B/{I^n})}{HH^l_{i}(B)}} \to
{\frac {HH^l_{i}(B/{I^{n_0}})}{HH^l_{i}(B)}}$ is zero for all $n \gg n_0$.
Now by comparing this exact sequence for $n_0$ and $n \ge n_0$ and 
using the doubling trick, we conclude the proof of ~\ref{eqn:AQH20}.

To prove ~\ref{eqn:AQH21}, we observe in ~\ref{eqn:AQH22} that
the right vertical arrow is injective for all $n \gg n_0$ by
Corollary~\ref{cor:AQH1} (part $(iii)$). Hence we can replace
the kernel of the middle vertical arrow by the kernel of the left
vertical arrow in order to prove  ~\ref{eqn:AQH21}, in which case
it follows from Corollary~\ref{cor:AQH1} (part $(i)$).
\end{proof}  
\begin{cor}\label{cor:DHC} 
Let the notations be as in  Corollary~\ref{cor:AQH1}. Then for any
conducting ideal $I$ for the normalization and for any $i \in \Z$ and
$n_0 \ge 0$, the natural maps of double relative Hochschild and cyclic
homology groups
\begin{equation}\label{eqn:DHC0} 
HH^l_i(A, B, I^n) \to HH^l_i(A, B, I^{n_0})
\end{equation}
\begin{equation}\label{eqn:DHC1} 
HC^l_i(A, B, I^n) \to HC^l_i(A, B, I^{n_0})
\end{equation}
are zero for all $n \gg n_0$.
\end{cor}
\begin{proof} The double relative Hochschild and cyclic homology are
classically known to be zero in negative degrees. So we assume 
$i \ge 0$. We first prove the result for the Hochschild homology.
The long exact double relative Hochschild homology sequence
gives for any $i, n \ge 0$, the short exact sequence
\[0 \to {\frac {HH^l_{i+1}(B, I^n)}{HH^l_{i+1}(A, I^n)}} \to
HH^l_i(A, B, I^n) \to {\rm Ker}\left(HH^l_i(A, I^{n}) \to 
HH^l_i(B, I^{n})\right) \to 0.\]
Comparing this exact sequence for $n_0$ and $n \ge n_0$, using 
Corollary~\ref{cor:AQH2} and the doubling trick, we get the desired
result.

To prove the result for the cyclic homology, we use induction on $i \ge 0$.
\\
We have $HC^l_i(A, B, I^n) = 0$ for $i < 0$ as pointed above, and 
hence the SBI-sequence give isomorphism
$HH^l_0(A, B, I^n) \xrightarrow {\cong} HC^l_0(A, B, I^n)$. So the
result holds for $i \le 0$. Suppose the result holds for all $j \le i-1$
with $i \ge 1$. We have the long exact SBI-sequence 
\[HH^l_i(A, B, I^n) \xrightarrow {I} HC^l_i(A, B, I^n) 
\xrightarrow {S} HC^l_{i-2}(A, B, I^n).\]   
Comparing this exact sequence for $n_0$ and $n \ge n_0$, using induction
on $i$ and ~\ref{eqn:DHC0} (which we just proved), 
the doubling trick gives us the proof of ~\ref{eqn:DHC1}.
\end{proof}

{\bf Proof of Theorem~\ref{thm:Artin-Rees}}:} Let $X$ be a quasi-projective
variety of dimension $d$ over a field $k$ and let $f : {\ov X} \to X$
be the smooth normalization of $X$. Let $Y \inj X$ be a given conducting
subscheme for the normalization. 

For $i \in \Z$, let ${\sK}_{i, (X, {\ov X}, Y)}$ denote the sheaf of double 
relative $K$-groups on $X$. This is a sheaf whose stalk of at any
point $x \in X$ is the double relative $K$-group 
$K_i({\sO}_{X, x}, {\sO}_{{\ov X}, x}, {\sI}_{Y,x})$, where ${\sI}_Y$ is
the ideal sheaf of $Y$.  By \cite{PW}, Proposition~A.5 (see also  
\cite{TT} for more general Brown-Gersten spectral sequences),
there exists a convergent spectral sequence
\[{^n}E^{p,q}_2 = H^p_{\rm Zar}\left(X, {\sK}_{q, (X, {\ov X}, nY)}\right) 
\Rightarrow K_{q-p}(X, {\ov X}, nY),\] with differential
$d_r : {^n}E^{p,q}_r \to {^n}E^{p+r,q+r-1}_r$.
This gives a finite filtration
\[K_i(X, {\ov X}, nY) = F^0_n \supseteq F^1_n \supseteq \cdots \supseteq 
F^d_n \supseteq F^{d+1}_n = 0\] such that for $0 \le j \le d$,
${F^j_n}/{ F^{j+1}_n} \cong {^n}E^{j,i+j}_{\infty} = {^n}E^{j,i+j}_{r}$,
where $r$ depends only on $d$ and $0 \le j \le d$. 
Thus by the descending induction on $j$ and using the doubling trick, it
suffices to show that for $0 \le j \le d$ and $n_0 \ge 0$, the natural map
${^n}E^{j,i+j}_{2} \to {^{n_0}}E^{j,i+j}_{2}$ is zero for all $n \gg n_0$.
This will be proved if we show that for $i \in \Z$ and $n_0 \ge 0$,
the natural map of sheaves ${\sK}_{i, (X, {\ov X}, nY)} \to
{\sK}_{i, (X, {\ov X}, {n_0}Y)}$ is zero for all $n \gg n_0$. For $ i \le 0$,
these sheaves are classically known to be zero (\cite{Bass}). So we assume
$i \ge 1$. Since $X$ is a $k$-variety, it is enough to show this when $X$
is affine. Thus we need to show that if $A$ is an essentially of
finite type $k$-algebra and $B$ is the smooth normalization of $A$, then for 
any conducting ideal $I$ and $i, n_0 \ge 0$, the natural map
$K_i(A, B, I^n) \to K_i(A, B, I^{n_0})$ is zero for all $n \gg n_0$.
But this follows immediately from Corollary~\ref{cor:DHC} (with
$l = \Q$) and Cortinas' result (\cite{Cortinas}, Corollary~0.2) that
these double relative $K$-groups are in fact rational vector spaces
and the Chern character map $K_i(A, B, I^n) \to HC^{\Q}_{i-1}(A, B, I^n)$
is an isomorphism. This completes the proof of Theorem~\ref{thm:Artin-Rees}.
$\hfill \square$ 
\\
\section{Formula for The Chow group of Zero Cycles}
Let $X$ be a quasi-projective variety of dimension $d \ge 2$ over a field $k$.
We assume in this section that $X$ is Cohen-Macaulay (all local rings
are Cohen-Macaulay) and it has only isolated singularities. Note that
this automatically implies that $X$ is normal. Our aim in this section
is to prove Theorem~\ref{thm:Chowgroup}. 
So let $p : {\wt X} \to X$ be a good resolution of singularities of $X$
and let $E$ denote the reduced exceptional divisor on $\wt X$. Let 
$S \subset X$ be the singular locus of $X$. We give $S$ the reduced induced
subscheme structure and denote by $nS$, the $n$th infinitesimal thickening
of $S$ in $X$. Let $Y \inj X$ be a closed subscheme of $X$ such that $p$
is the blow-up of $X$ along $Y$. Then $S = Y_{red}$. Let $\sI$ denote
the ideal sheaf for $Y$. Then one has
\[\wt {X} = {\rm Proj}_{X}( \oplus_{n \ge 0} {\sI}^n)\ \ {\rm and}\ \ 
E = ({\rm Proj}_{Y}( \oplus _{n \ge 0} {\sI}^n / {\sI}^{n+1}))_{red}.\]
Putting ${\wt Y} = {\rm Proj}_{Y}( \oplus _{n \ge 0} 
{\sI}^n /{\sI}^{n+1})$, we get $S \subset Y \subset nS$ and $ E \subset 
{\wt Y} \subset nE$ for all sufficiently large $n$.
\begin{lem}\label{lem:commute} 
For all $n \ge 1$, the map $F^dK_0(X, nS) \to F^dK_0(X)$ is
an isomorphism. In particular, there are natural maps 
$F^dK_0(X) \to F^dK_0({\wt X}, nE)$ such that Diagram ~\ref{eqn:relative}
commutes and all maps there are surjective. 
\end{lem}
\begin{proof} Since $S$ is zero-dimensional, the map 
$F^dK_0(X, nS) \to F^dK_0(X)$ is an isomorphism by \cite{Kr1}
(Lemma~3.1). On the other hand,  there are natural maps
$F^dK_0(X, nS) \to F^dK_0({\wt X}, nE)$  by the definition 
of relative $K$-groups (see Section~1).
The isomorphism above now shows that this map
factors through a map $F^dK_0(X) \to F^dK_0({\wt X}, nE)$.
The surjectivity assertion follows from Lemmas ~3.1 and ~3.2 of \cite{Kr1}.
\end{proof}

{\bf Proof of Theorem~\ref{thm:Chowgroup}}:} The Northcott-Rees theory 
gives a minimal reduction of ideal sheaf 
${\sJ} \subset {\sI}$ of ${\sI}$ in the sense that 
${\sJ}{{\sI}^n} = {\sI}^{n+1}$ for all sufficiently large $n$.
Furthermore, since $X$ is Cohen-Macaulay and $S$ is a finite set of closed
points, we can choose ({\sl cf.} \cite{Weibel2}) ${\sJ}$ to be a local 
complete intersection ideal sheaf on $X$.  
Now we follow the proof of Theorem~1.1 of \cite{KS}
to get a factorization
\[
\xymatrix{
\wt {X} \ar [drr] ^{f} \ar [dd] _{p} \\
&& X' \ar [dll] ^{p'}\\
X && }
\]
where $p'$ is the blow-up of $X$ along ${\sJ}$ and $f$ is the
normalization morphism. Let $Y_1$ denote the local complete intersection 
subscheme of $X$ defined by ${\sJ}$. Since ${\sJ}$ is a reduction for 
${\sI}$, we see that
${Y_1}_{red} = Y_{red} = S$ and hence ${\sI}^n \subset {\sJ} \subset
{\sI}$ for all large $n$. Let 
\[Y' = Y \times_X X', \ \ \ \ \  Y_1' = Y_1 \times_X X', \ \
\ \ \ {\wt Y} = Y \times_X {\wt X},\]
\[{\wt Y}_1 = Y_1 \times_X {\wt X},
\ \ \ \ \ \rm{and} \ \ S' = {(S \times_X X')}_{red}.\]
Let $Z' \subset X'$ be a conducting subscheme for the normalization map
$f$. Put ${\wt Z}' = Z' {{\times}_{X'}} {\wt X}$. Then we see that 
${Z'}_{red} \subset S'$ and ${\wt Z}'_{red} \subset E$. In particular, 
given any $m > 0$,
we have $mZ' \subset nS'$ and $m{\wt Z}' \subset nE$ for all large $n$.
Hence for a given $m > 0$, we have the following commutative diagram for all 
sufficiently large $n$ with all maps surjective.

\begin{equation}\label{eqn:relative0}
\xymatrix{
F^dK_0( \wt {X},nE) \ar @{->>}[rr] &&
F^dK_0(\wt {X}, m{\wt Z}') \ar @{->>}[rr] &&
F^dK_0(\wt {X}) \\
F^dK_0(X',nS') \ar @{->>}[rr] \ar @{->>}[u] &&
F^dK_0(X',mZ') \ar @{->>}[rr] \ar @{->>}[u] &&
F^dK_0(X{'}) \ar @{->>}[u] \\
F^dK_0(X,nS) \ar @{->>}[u] \ar @{->>}[rr] &&
F^dK_0(X) \ar @{->>}[urr] \ar @{->>}[ull] && 
}
\end{equation}
The surjectivity of all maps follows from Lemma~3.2 of \cite{Kr1}, 
and the bottom horizontal map is an isomorphism by Lemma~\ref{lem:commute}.
Now since $p'$ is a blow-up along a local complete intersection subscheme,
the map $F^dK_0(X) \to F^dK_0(X')$ is also injective by 
\cite{KS} (Corollary 2.5). Combining this
with the surjectivity of arrows in the above diagram, we get another 
diagram below with all the arrows being isomorphisms.
\[
\xymatrix{
F^dK_0(X',nS') \ar[rr] &&
F^dK_0(X',mZ') \ar[rr] &&
F^dK_0(X{'}) \\
F^dK_0(X,nS) \ar[u] \ar[rr] &&
F^dK_0(X) \ar[urr] \ar[ull] && 
}
\]
Next we study the relation between $F^dK_0(X',mZ')$ and 
$F^dK_0(\wt {X}, m{\wt Z'})$ for fixed $m > 0$. By the long exact sequence
of double relative $K$-theory, one has an exact sequence
\[K_0(X', {\wt X}, mZ') \to K_0(X', mZ') \to K_0({\wt X}, m{\wt Z}').\]  
We compare this exact sequence for $m = 1$ and $m >> 0$ to get a diagram
\[
\xymatrix@C.6pc{
K_0(X', {\wt X}, mZ') \ar[r] \ar[d] &  K_0(X', mZ') \ar[r] \ar[d] &
K_0({\wt X}, m{\wt Z}') \ar[d] \\
K_0(X', {\wt X}, Z') \ar[r] &  K_0(X', Z') \ar[r] &
K_0({\wt X}, {\wt Z}')
}
\]
The left vertical map is zero for $m \gg 0$ by 
Theorem~\ref{thm:Artin-Rees}. \\
Put $A_m = {\rm Ker}\left(F^dK_0(X',mZ') \to 
F^dK_0({\wt X}, m{\wt Z}')\right)$. 
Then the above diagram gives another diagram of short exact sequences
\[
\xymatrix@C.6pc{
0 \ar [r] & A_m \ar [r] \ar [d] & F^dK_0(X',mZ') \ar [r] \ar [d] &
F^dK_0({\wt X}, m{\wt Z}') \ar [r] \ar [d] & 0 \\
0 \ar [r] & A_1 \ar [r] &
F^dK_0(X',Z') \ar [r] & F^dK_0({\wt X}, {\wt Z}') \ar [r] & 0
}
\]
where the left vertical map is zero for $m \gg 0$. 
By mapping Diagram ~\ref{eqn:relative0} to a similar
diagram with $m = 1$, we see that the middle
vertical map above is an isomorphism. A diagram chase
shows that $F^dK_0(X',mZ') \to F^dK_0({\wt X}, m{\wt Z}')$ is an isomorphism.
This, together with the isomorphism $F^dK_0(X',nS') \rightarrow 
F^dK_0(X',mZ')$ now shows that the map $F^dK_0(X,nS) \to    
F^dK_0({\wt X}, nE)$ is an isomorphism for all large $n$. 
This gives the desired isomorphisms
\[F^dK_0({\wt X}, nE) \to F^dK_0({\wt X}, {(n-1)}E) \ {\rm and} \
F^dK_0(X) \to F^dK_0(\wt {X}, nE)\] for all large $n$. 
Finally, for $X$ affine or projective, 
we get 
\[CH^d(X) \cong~ ~{\underset {n}{\varprojlim}}~F^{d}K_{0}({\wt X},nE)\] 
if $k$ is algebraically closed.
$\hfill \square$
\section{Cohomology of Milnor $K$-sheaves}
It is by now a well known fact that algebraic cycles are closely connected
to the cohomology of Quillen $K$-theory sheaves. This is true even for
certain classes of singular varieties. However, the Quillen $K$-theory groups
are often very difficult to compute. On the other
hand, one also has the sheaves of Milnor $K$-theory on schemes which are
relatively simpler looking objects. Our applications of 
Theorem~\ref{thm:Chowgroup} in this paper are based on the observation
that for the purposes of zero cycles, the appropriate cohomology of 
Quillen $K$-theory sheaves can often be approximated by the cohomology of 
Milnor $K$-theory sheaves, which can be computed by 
some other means. In this section, we prove certain general reduction
steps in order to use Theorem~\ref{thm:Chowgroup} for 
studying the Chow group of zero cycles on varieties with isolated
singularities. In particular, we give a very precise sufficient condition
for the Chow group of zero cycles on a variety $X$ with Cohen-Macaulay
isolated singularities to be isomorphic to the similar group on a
resolution of singularities of $X$. We begin with the following result.

For any variety $X$ over a field $k$, let ${\sK}^M_{m, X}$
denote the sheaf of Milnor $K$-groups on $X$. This is a sheaf
whose stalk at any point $x$ of $X$ is the Milnor $K$-group of the local
ring ${\sO}_{X, x}$. For any closed embedding $i : Y \inj X$, let 
${\sK}^M_{m, (X, Y)}$ 
be the sheaf of relative Milnor $K$-groups defined so that the sequence
of sheaves
\begin{equation}\label{eqn:M1} 
0 \to {\sK}^M_{m, (X, Y)} \to {\sK}^M_{m, X} \to
i_*({\sK}^M_{m, Y}) \to 0
\end{equation}
is exact. Note that the map
${\sK}^M_{m, X} \to i_*({\sK}^M_{m, Y})$ is always surjective.
From now on, we shall assume our base field $k$ to algebraically closed
unless mentioned otherwise.
\begin{lem}\label{lem:iso}    
Let $X$ be an affine or projective variety of dimension $d$ over $k$. Then
there are natural isomorphisms
\[CH^d(X) \cong F^dK_0(X) \cong H^d(X, {\sK}_{d, X}) \cong  
H^d(X, {\sK}^M_{d, X}).\]
\end{lem}
\begin{proof} The isomorphism $CH^d(X) \cong F^dK_0(X)$ was shown by
Levine (\cite{Levine1}, Corollary~2.7 and Theorem~3.2). In this case,
Barvieri-Viale has shown (\cite{Viale}, Corollary~A) that
there is a natural surjection $CH^d(X) \surj H^d(X, {\sK}_{d, X})$
with finite kernel. If $X$ is affine, then this kernel must be zero
by \cite{Levine1} (Theorem~2.6). If $X$ is projective, then
there is an albanese map
$CH^d(X) \to H^d(X, {\sK}_{d, X}) \to H^d({\wt X}, {\sK}_{d, {\wt X}})
\to Alb(\wt X)$, where $\wt X$ is a resolution of singularities of $X$. 
Now the isomorphism $Alb(X) {\cong} Alb(\wt X)$ (since $X$ is normal)
and the Roitman torsion theorem
implies that $CH^d(X) \to H^d(X, {\sK}_{d, X})$ must be an isomorphism.
Finally, the isomorphism $H^d(X, {\sK}_{d, X}) \cong  
H^d(X, {\sK}^M_{d, X})$ is shown in \cite{Kr1} (Corollary~4.2).
\end{proof}
\begin{prop}\label{prop:Milnor}
Let $\wt X$ be a smooth quasi-projective variety of dimension $d+1$ over
a field $k$. Let $E \inj {\wt X}$ be a strict normal crossing divisor.
Then the natural cup product map
\[H^d(E, {\sK}^M_{m, E}) \otimes k^* \to H^d(E, {\sK}^M_{{m+1}, E})\]
is surjective for all $m \ge d$.
\end{prop}
%$ m \ge 1$? 
\begin{proof} We prove this by induction on $d \ge 1$ and divide the
proof into several cases. \\
{\bf Case I:} $d$ arbitrary and $E$ smooth. \\
In this case one has the Gersten resolution
(\cite{MSV}, Proposition~4.3)
\[{\sK}^M_{m, E} \to i_*\left(K^M_m (k(E))\right) \to
{\underset {x \in E^{(1)}} {\oplus}} i_*\left(K^M_{m-1}(k(x))\right)
\to \cdots \hspace*{3cm}\]
\[\hspace*{3cm} \to 
{\underset {x \in E^{(d-1)}} {\oplus}} i_*\left(K^M_{m-d+1}(k(x))\right)
\to 
{\underset {x \in E^{(d)}} {\oplus}} i_*\left(K^M_{m-d}(k(x))\right)
\to 0, \]
where the first map is generically injective (in fact everywhere
injective by a recent result of Kerz \cite{Kerz}). This resolution gives
a commutative diagram
\[
\xymatrix@C.5pc{
\left({\underset {x \in E^{(d)}} {\oplus}} i_*\left(K^M_{m-d}(k(x))\right)
\right) {\otimes} k^{*} \ar@{->>}[r] \ar[d] &
H^d(E, {\sK}^M_{m, E}) \otimes k^* \ar[d] \\
{\underset {x \in E^{(d)}} {\oplus}} i_*\left(K^M_{m+1-d}(k(x))\right)
\ar@{->>}[r] & H^d(E, {\sK}^M_{{m+1}, E}).}
\]
The horizontal arrows in this diagram are surjective by the above
resolution. The left vertical arrow is surjective since $k(x) = k$
for $x \in E^{(d)}$ as $k$ is algebraically closed and the map
$K^M_i(k) \otimes k^M_j(k) \to K^M_{i+j} (k)$ is surjective. A diagram 
chase proves the result. \\
\\
{\bf Case II:} $d = 1$ and $E$ not smooth. \\
Let $E = E_1 \cup \cdots \cup E_r$ with $r \ge 2$ and put
$E' = E_1 \cup \cdots \cup E_{r-1}$, $F_r = E' \cap E_r$. Since $E$ is
a strict normal crossing divisor on $\wt X$, we see that $F_r$ is a
strict normal crossing divisor on $E_r$. Thus $F_r$ is finite set
of closed points for which the proposition is obvious.
Let ${\sI}_{E'}$ (resp ${\sI}_{E_r}$) denote the ideal sheaf of
$E'$ (resp $E_r$) on $\wt X$. Let ${\ov E}$ be the closed subscheme of
$\wt X$ defined by the sheaf of ideals ${\sI}_{E'} \cap {\sI}_{E_r}$.
Then one has 
\begin{equation}\label{eqn:Milnor0}
{\sK}^M_{m, E} \surj  {\sK}^M_{m, {\ov E}}
\end{equation}
and this map is generically an isomorphism. 
In particular, we have 
\begin{equation}\label{eqn:Milnor01}
H^d(E, {\sK}^M_{m, E}) \xrightarrow {\cong}
H^d(E, {\sK}^M_{m, {\ov E}}) \ \forall \ m.
\end{equation}
Note here that
$E$ is the subscheme of $\wt X$ locally defined by the product of
${\sI}_{E'}$ and ${\sI}_{E_r}$.   
By \cite{Kr1} (Lemma~4.5), there is a short exact sequence of sheaves
\begin{equation}\label{eqn:Milnor1}
0 \to {\sK}^M_{m, {\ov E}} \to  i_*\left({\sK}^M_{m, E'}\right)
\oplus i_*\left({\sK}^M_{m, E_r}\right) \to 
i_*\left({\sK}^M_{m, F_r}\right) \to 0. 
\end{equation}
Taking the cohomology exact sequences, we get a commutative diagram
of exact sequences

\begin{equation}\label{eqn:Milnor2}
\xymatrix@R=1.5cm{
H^{d-1}(F_r, {\sK}^M_{m, F_r}) \otimes k^* \ar[d] \ar[r] &
H^d(E, {\sK}^M_{m, {\ov E}}) \otimes k^* \ar[d] \ar[r] &
\save[]+<0cm,.5cm>*\txt{$H^d(E', {\sK}^M_{m, E'}) \otimes k^*$}\restore
\text{\ \ \ \ \ \ \ \ \ \ \ \ }\oplus\text{\ \ \ \ \ \ \ \ \ \ \ \ \ }
 \save[]-<0cm,.5cm>="src"*\txt{$H^d(E_r, {\sK}^M_{m, E_r}) 
\otimes k^*$}\restore
 \ar[r] & 0 \\
H^{d-1}(F_r, {\sK}^M_{m+1, F_r}) \ar[r] &
H^d(E, {\sK}^M_{m+1, {\ov E}}) \ar[r] &
\save[]+<0cm,.5cm>="tar"*\txt{$H^d(E', {\sK}^M_{m+1, E'})$}\restore 
\text{\ \ \ \ \ \ \ \ \ \ \ \ }\oplus\text{\ \ \ \ \ \ \ \ \ \ \ \ \ }
\save[]-<0cm,.5cm>*\txt{$H^d(E_r, {\sK}^M_{m+1, E_r})$}\restore
\ar[r] & 0,
\ar "src"-<0cm,.275cm>;"tar"+<0cm,.275cm>
}
\end{equation}
where the last terms in both rows are zero because dim$(F_r) \le d-1$.
We have already observed that the left vertical arrow is surjective.
The map $H^d(E_r, {\sK}^M_{m, E_r}) \otimes k^* \to  
H^d(E_r, {\sK}^M_{m+1, E_r})$ is surjective by Case I. 
The map between the other summand of the right vertical arrow is surjective
by induction on the number of components since $n(E') < n(E)$. 
We complete the proof of Case II by a diagram chase and ~\ref{eqn:Milnor01}. \\
\\
{\bf Case III:} $d \ge 2$.
Assume by induction that the proposition holds for whenever dimension
of the normal crossing divisor is $d' < d$. Let $n(E)$ denote the
number of irreducible components of $E$. We now induct on $n(E)$.
The case $n(E) = 1$ is already proved above. So assume $n(E) = n \ge 2$
and assume we have proved Case II for all normal crossing divisors
$E'$ with $n(E') < n$. 
Let $E', E_r, F_r$ and $\ov E$ be as in Case II. Then $E'$ is a strict
normal crossing divisor on $\wt X$ with $n(E') < n$. Moreover,
$E_r$ is smooth and $F_r$ is a strict normal crossing divisor on $E_r$
(if not empty) and $E_r$ is a smooth variety of smaller dimension than 
$\wt X$. 
So the proposition holds for $E', F_r$ and $E_r$ by induction and smooth
case.
A diagram chase again in ~\ref{eqn:Milnor2} and ~\ref{eqn:Milnor01} now 
complete the proof.
\end{proof}
\begin{prop}\label{prop:Milnor3}
Let $Z$ be a quasi-projective variety of dimension $d$ over a field $k$ 
and let $W= Z_{\rm red}$. For $i \ge 0$, let ${\Omega}^i_{(Z, W)/{\Q}}$
denote the kernel of the map ${\Omega}^i_{Z/{\Q}} \to 
{\Omega}^i_{W/{\Q}}$. Then there is a natural isomorphism
\[H^d\left(Z, {\sK}^M_{{d+1}, (Z, W)}\right) \to 
H^d\left(Z, {\frac {{\Omega}^d_{(Z, W)/{\Q}}} 
{{\del}\left({\Omega}^{d-1}_{(Z, W)/{\Q}}\right)}}\right),\]
where ${\del} : {\Omega}^i_{Z/{\Q}} \to {\Omega}^{i+1}_{Z/{\Q}}$
is the differential map.
\end{prop}
\begin{proof} Let ${\phi}_Z : {\sK}^M_{m, Z} \to {\sK}_{m, Z}$ denote the 
natural map from Milnor $K$-theory to Quillen $K$-theory. Since we are in 
characteristic 0, there exists by \cite{ST}(Theorem~12.3) a natural map 
\[{\psi}_Z : {\sK}_{m, Z} \to {\sK}^M_{m, Z}\]
such that ${\psi}_Z \circ {\phi}_Z = {({(-1)}^{m-1} {(m-1)}!)}{\rm Id} \ 
\forall \ m \ge 1$.
Furthermore, with respect to the $\gamma$-filtration and Adams operation
on the Quillen $K$-theory (\cite{Levine2}), one has
$F^{d+1}{\sK}_{d+1, Z} = {\sK}^{(d+1)}_{d+1, Z}$ and the map
${\phi}_Z : {\sK}^M_{d+1, Z} \to {\sK}^{(d+1)}_{d+1, Z}$ is isomorphism
modulo $d!$. In particular, we get natural maps 
\begin{equation}\label{eqn:Milnor30}
{\sK}^M_{d+1, Z} \xrightarrow {{\phi}_Z} {\sK}^{(d+1)}_{d+1, Z}
\xrightarrow {{\psi}_Z} {\sK}^M_{d+1, Z} 
\end{equation}
which are isomorphisms modulo $d!$.   

We now consider the commutative diagram
\[
\xymatrix@C.8pc{
0 \ar[r] &  {{\sK}^M_{d+1, (Z, W)}} \ar[r] \ar[d]^{{\phi}_{ZW}} &
{{\sK}^M_{d+1, Z}} \ar[d]^{{\phi}_Z} \ar[r] & 
{{\sK}^M_{d+1, W}} \ar[d]^{{\phi}_W} \ar[r] & 0 \\
0 \ar[r] &  {\ov {{\sK}^{(d+1)}_{d+1, (Z, W)}}} \ar[r] \ar[d]^{{\psi}_{ZW}} &
{{\sK}^{(d+1)}_{d+1, Z}} \ar[r] \ar[d]^{{\psi}_Z} & 
{\ov {{\sK}^{(d+1)}_{d+1, W}}} \ar[r] \ar[d]^{{\psi}_W} & 0 \\
0 \ar[r] &  {{\sK}^M_{d+1, (Z, W)}} \ar[r] &
{{\sK}^M_{d+1, Z}} \ar[r] & 
{{\sK}^M_{d+1, W}} \ar[r] & 0, }
\]
where the top and the bottom rows are exact by ~\ref{eqn:M1} and
the group ${\ov {{\sK}^{(d+1)}_{d+1, (Z, W)}}}$ is defined to make 
the middle row exact. A diagram chase gives maps ${\psi}_{ZW}$ and 
${\phi}_{ZW}$ which are inverses to each other modulo $d!$.
Since the map ${\sK}^{(d+1)}_{d+1, (Z, W)} \surj 
{\ov {{\sK}^{(d+1)}_{d+1, (Z, W)}}}$
is isomorphism on the smooth locus of $W$, we get 
$H^d\left(Z, {\sK}^{(d+1)}_{d+1, (Z, W)}\right) \cong
H^d\left(Z, {\ov {{\sK}^{(d+1)}_{d+1, (Z, W)}}}\right)$, which in turn shows 
that the map
\begin{equation}\label{eqn:Milnor31} 
H^d\left(Z, {\sK}^{(d+1)}_{d+1, (Z, W)}\right) \to
H^d\left(Z, {\sK}^M_{d+1, (Z, W)}\right)
\end{equation}
is an isomorphism modulo $d!$.
However, $W$ being the reduced part of $Z$,  
${\sK}^{(d+1)}_{d+1, (Z, W)}$ is a sheaf of $\Q$-vector spaces 
(\cite{Cath}, Section~1) and ${\sK}^M_{d+1, (Z, W)}$ is a sheaf of
divisible groups by \cite{Kr1} (Sublemma~4.8). In particular, the
left hand side in ~\ref{eqn:Milnor31} is a $\Q$-vector space and the
right hand side is a divisible group (since $H^d$ is right exact on 
$Z$). This shows that the map in ~\ref{eqn:Milnor31} must be an
isomorphism.

By \cite{Cath} (Theorem~1), there exists a functorial isomorphism of
filtered sheaves of $\Q$-vector spaces 
${\sK}_{d+1, (Z, W)} \xrightarrow {\cong} {\sH}{\sC}_{d, (Z, W)}$,
where the filtration is given by the $\gamma$-filtration on both sides,
and ${\sH}{\sC}$ denote the sheaves of cyclic homology over
the base $\Q$.
This gives an isomorphism 
${\sK}^{(d+1)}_{d+1, (Z, W)} \xrightarrow {\cong} 
{\sH}{\sC}^{(d)}_{d, (Z, W)}$ and by ~\ref{eqn:Milnor31}, we get
$H^d\left(Z, {\sK}^M_{d+1, (Z, W)}\right) \cong
H^d\left(Z, {\sH}{\sC}^{(d)}_{d, (Z, W)}\right)$. Thus we are reduced
to showing that there is a natural isomorphism
\begin{equation}\label{eqn:Milnor32}
H^d\left(Z, {\sH}{\sC}^{(d)}_{d, (Z, W)}\right)
\xrightarrow {\cong} 
H^d\left(Z, {\frac {{\Omega}^d_{(Z, W)/{\Q}}} 
{{\del}\left({\Omega}^{d-1}_{(Z, W)/{\Q}}\right)}}\right).
\end{equation}
We have an exact sequence of sheaves
\[{\sH}{\sC}^{(d)}_{d, (Z, W)} \to {\sH}{\sC}^{(d)}_{d, Z} \to
{\sH}{\sC}^{(d)}_{d, W} \to 0,\]
where the last term is zero since
\[{\sH}{\sC}^{(d)}_{d, Z} \cong {\frac {{\Omega}^d_{Z}/{\Q}}  
{{\del} \left({\Omega}^{d-1}_{Z}/{\Q}\right)}} \surj
{\frac {{\Omega}^d_{W}/{\Q}} {{\del} \left({\Omega}^{d-1}_{W}/{\Q}\right)}}
\cong {\sH}{\sC}^{(d)}_{d, W}\]
by \cite{Loday} (Theorem~4.6.8). Furthermore, since ${\sO}_Z
\surj {\sO}_W$ is locally split on the smooth locus of $W$, we
see that the first map in the above exact sequence is injective on
the smooth locus of $W$. By the same reason, there is a natural
surjection 
\[{\frac {{\Omega}^d_{(Z, W)/{\Q}}} 
{{\del}\left({\Omega}^{d-1}_{(Z, W)/{\Q}}\right)}} \surj
{\rm Ker}\left({\frac {{\Omega}^d_{Z}/{\Q}}  
{{\del} \left({\Omega}^{d-1}_{Z}/{\Q}\right)}} \surj
{\frac {{\Omega}^d_{W}/{\Q}} 
{{\del} \left({\Omega}^{d-1}_{W}/{\Q}\right)}}\right),\]
which is an isomorphism on the smooth locus of $W$. In particular,
we get surjective maps
\[
\xymatrix@C.9pc{
{\sH}{\sC}^{(d)}_{d, (Z, W)} \ar@{->>}[r] &
{\rm Ker}\left({\frac {{\Omega}^d_{Z}/{\Q}}  
{{\del} \left({\Omega}^{d-1}_{Z}/{\Q}\right)}} \surj
{\frac {{\Omega}^d_{W}/{\Q}} {{\del} \left({\Omega}^{d-1}_{W}/{\Q}\right)}}
\right) &
{\frac {{\Omega}^d_{(Z, W)/{\Q}}} 
{{\del}\left({\Omega}^{d-1}_{(Z, W)/{\Q}}\right)}} \ar@{->>}[l]}
\]
which are isomorphisms on the smooth locus of $W$, and hence they induce
isomorphisms on the top cohomology $H^d$.
This proves ~\ref{eqn:Milnor32} and hence the proposition.
\end{proof}
\begin{cor}\label{cor:sufficient}
Let $X$ be either an affine or a projective variety of dimension $d \ge 2$
over a field $k$. Assume that $X$ is Cohen-Macaulay and has only isolated 
singularities. Let $p : {\wt X} \to X$ be a good resolution of 
singularities of $X$ such that
the exceptional divisor $E$ is a strict normal crossings divisor 
(such a resolution always exists). Then the map $CH^d(X) \to CH^d(\wt X)$
is an isomorphism if the following two conditions hold. \\    
$(i)$ The map $H^{d-1}\left({\wt X}, {\sK}_{{d-1}, {\wt X}}\right) 
{\otimes} k^*
\to H^{d-1}\left(E, {\sK}_{{d-1}, E}\right) {\otimes} k^*$ is surjective. \\
$(ii)$ $H^{d-1}\left(nE, {\frac {{\Omega}^{d-1}_{(nE, E)/{\Q}}} 
{{\del}\left({\Omega}^{d-2}_{(nE, E)/{\Q}}\right)}}\right) = 0$
for all $n \ge 1$.
\end{cor}
\begin{proof} By Theorem~\ref{thm:Artin-Rees}, we need to show that
the map $F^dK_0({\wt X}, nE) \to F^dK_0(\wt X)$ is an isomorphism
for all large $n$. By \cite{Kr1} (Proposition~4.3), this further
reduces to showing that the map 
$H^d\left({\wt X}, {\sK}^M_{d, ({\wt X}, nE)}\right) \to
H^d\left({\wt X}, {\sK}^M_{d, {\wt X}}\right)$ is an isomorphism
for all $n \ge 1$.   
Considering the long exact cohomology sequence corresponding to ~\ref{eqn:M1},
\[H^{d-1}\left({\wt X}, {\sK}^M_{d, {\wt X}}\right) \to 
H^{d-1}\left(nE, {\sK}^M_{d, nE}\right) \to
H^d\left({\wt X}, {\sK}^M_{d, ({\wt X}, nE)}\right) \to 
H^d\left({\wt X}, {\sK}^M_{d, {\wt X}}\right) \to  0,\]
we need to show that the first map on the left is surjective for all $n$.

First we consider the case $n =1$. We have seen in the proof of
Proposition~\ref{prop:Milnor3} that there is a natural map
${\psi}_E : {\sK}_{{d-1}, E} \to {\sK}^M_{{d-1}, E}$ whose cokernel is
of fixed exponent $(d-2)!$. In particular, the cokernel of the
map $H^{d-1}\left(E, {\sK}_{{d-1}, E}\right) \to
H^{d-1}\left(E, {\sK}^M_{{d-1}, E}\right)$ is of finite exponent
(since $H^{d-1}$ is right exact on $E$). However, as $k^*$ is a divisible
group ($k = {\ov k}$), we must have
\begin{equation}\label{eqn:suff0}
H^{d-1}\left(E, {\sK}_{{d-1}, E}\right) {\otimes} k^*
\surj H^{d-1}\left(E, {\sK}^M_{{d-1}, E}\right) {\otimes} k^*.
\end{equation}

Next we have a commutative diagram
\[
\xymatrix{
H^{d-1}\left({\wt X}, {\sK}^M_{{d-1}, {\wt X}}\right) {\otimes} k^* \ar[d]
\ar[r] & H^{d-1}\left({\wt X}, {\sK}^M_{d, {\wt X}}\right) \ar[d] \\
H^{d-1}\left(E, {\sK}^M_{{d-1}, E}\right) {\otimes} k^*
\ar[r] & H^{d-1}\left(E, {\sK}^M_{d, E}\right).}
\]
The first condition of the corollary and ~\ref{eqn:suff0} together
imply that the left vertical arrow
is surjective. The bottom horizontal arrow is surjective by
Proposition~\ref{prop:Milnor}. A diagram chase shows that the right vertical
arrow is also surjective. This finishes the case $n = 1$.

Now assume $n \ge 2$. Then ~\ref{eqn:M1} for the pair $E \inj nE$
gives an exact sequence 
\[H^{d-1}\left(nE, {\sK}^M_{d, (nE, E)}\right) \to 
H^{d-1}\left(nE, {\sK}^M_{d, nE}\right) \to 
H^{d-1}\left(E, {\sK}^M_{d, E}\right) \to  0.\]
The second condition of the corollary and Proposition~\ref{prop:Milnor3}
together imply that the group on the left is zero. This
completes the proof.
\end{proof}
\section{Chow group of affine cones}
Let $Y \inj {\P}^N_{k}$ be a smooth projective variety of dimension $d$
over a field $k$. 
Let $X = C(Y)$ be the affine cone over $Y$ and let ${\ov X} \inj
{\P}^{N+1}_k$ be the projective cone over $Y$. Let $P \in X \inj {\ov X}$
denote the vertex of the cone.
We assume that $X$ is Cohen-Macaulay. Since $P$ is the only singular
point of both $X$ and $\ov X$, we see that $\ov X$ is also Cohen-Macaulay. Let
$p: {\wt X} \to X$ be the blow-up of $X$ along the vertex $P$ and
let $E = p^{-1}(P)$ be the exceptional divisor for the blow-up.
This situation gives rise to the following commutative diagram.
\begin{equation}\label{eqn:blowup}
\xymatrix@C.6pc{   
E \ar@{^{(}->}[r]^{i} \ar[ddrr] & {\wt X} \ar[rr]^{p} \ar@{^{(}->}[dr] 
\ar[ddr]_{\pi} & & X \ar@{^{(}->}[dr] & \\
& & Z \ar[d]^{\ov {\pi}} \ar[rr]^{\ov p} & & {\ov X} \\
& & Y & &}
\end{equation} 
The map $\ov p$ is the blow-up of $\ov X$ along $P$. The map $\pi$
is the ${\A}^1$-bundle over $Y$ associated to the
ample line bundle ${\sO}_Y(1)$ with the zero-section $E$ and 
$\ov {\pi}$ is the ${\P}^1$-bundle over $Y$ associated to the vector bundle 
${\sO}_Y \oplus {\sO}_Y(1)$.
In particular, $\wt X$ (resp $Z$) is a good resolution of singularities 
of $X$ (resp $\ov X$) such that the exceptional divisor 
$E \inj {\wt X} \inj Z$ is smooth and the inclusion of $E$ in $\wt X$ and
$Z$ has sections given by the maps $\pi$ and $\ov {\pi}$. 
Our aim in this section is to prove Theorem~\ref{thm:cone} for which
we need the following preliminary results.
\begin{lem}\label{lem:prelim0}
Let $\wt X$ be a smooth quasi-projective variety of dimension $d+1$ over
a field $k$ and let $E \inj \wt X$ be a smooth divisor such that
\[H^d\left(E, {\Omega}^i_{E/k} {\otimes}_E {\frac {{\sI}^n}
{{\sI}^{n+1}}}\right)
= 0 \ {\rm for} \ i \ge 0 \ {\rm and} \ n \ge 1,\]
where $\sI$ is the ideal sheaf of $E$ on $\wt X$. Then
\[H^d\left(nE, {\Omega}^i_{(nE, E)/k}\right) = 0 
 \ {\rm for} \ i \ge 0 \ {\rm and} \ n \ge 1.\]
\end{lem}
\begin{proof} We first claim that
\begin{equation}\label{eqn:prelim01} 
H^d\left(nE, {\Omega}^i_{{nE}/k} {\otimes} {\frac {\sI}
{{\sI}^n}}\right) = 0 \ \forall \ i \ge 0, n \ge 1.
\end{equation}
To prove the claim, we consider for $i \ge 0$ and $n \ge 1$
the compatible maps
\[0 = {\Omega}^i_{{nE}/k} {\otimes} {\frac {{\sI}^n}
{{\sI}^n}} \to {\Omega}^i_{{nE}/k} {\otimes} {\frac {{\sI}^{n-1}}
{{\sI}^n}} \to \cdots \to  {\Omega}^i_{{nE}/k} {\otimes} {\frac {\sI}
{{\sI}^n}}\] and put
$F^i_j = {\rm Image} \left({\Omega}^i_{{nE}/k} {\otimes} {\frac {{\sI}^j}
{{\sI}^n}} \to {\Omega}^i_{{nE}/k} {\otimes} {\frac {\sI}
{{\sI}^n}}\right)$ for $1 \le j \le n$. This gives a finite
filtration ${\{F^i_j\}}_{1 \le j \le n}$ of 
${\Omega}^i_{{nE}/k} {\otimes} {\frac {\sI} {{\sI}^n}}$ such that
\begin{equation}\label{eqn:prelim02} 
{\Omega}^i_{nE/k} {\otimes}_{nE} {\frac {{\sI}^j}
{{\sI}^{j+1}}} \surj {\frac {F^i_j} {F^i_{j+1}}} \ \forall \ 1 \le j \le n.
\end{equation}
Since $H^d$ is right exact on $nE$, this filtration also gives exact sequence
\[H^d\left(nE, F^n_{j+1}\right) \to H^d\left(nE, F^n_{j}\right) \to
H^d\left(nE,{\frac {F^i_j} {F^i_{j+1}}}\right) \to 0\]
for $1 \le j \le n$. Thus by a descending induction on $j$ and by
~\ref{eqn:prelim02}, it suffices to show that
\begin{equation}\label{eqn:prelim03}
H^d\left(nE, {\Omega}^i_{nE/k} {\otimes}_{nE} {\frac {{\sI}^j}
{{\sI}^{j+1}}}\right) = 0  \ \forall \ 1 \le j \le n
\end{equation}
in order to prove the claim.  
If $n=1$ or $i = 0$, this is our assumption. So assume $n \ge 2$ and $i \ge 1$.
Then we have 
\[{\Omega}^i_{nE/k} {\otimes}_{nE} {\frac {{\sI}^j}
{{\sI}^{j+1}}} \cong {\Omega}^i_{nE/k} {\otimes}_{nE} \left({\sO}_E
{\otimes}_E {\frac {{\sI}^j} {{\sI}^{j+1}}}\right) \cong
{\ov {{\Omega}^i_{nE/k}}} {\otimes}_E {\frac {{\sI}^j} {{\sI}^{j+1}}},\]
where ${\ov {{\Omega}^i_{nE/k}}} = {\Omega}^i_{nE/k} {\otimes}_{nE} {\sO}_E$.
Next we have short exact sequence
\begin{equation}\label{eqn:prelim04} 
0 \to {\frac {\sI} {{\sI}^2}} \to {\ov {{\Omega}^1_{nE/k}}} \to 
{\Omega}^1_{E/k} \to 0,
\end{equation}
where the first term is zero because the long exact sequence of 
Andr{\'e}-Quillen homology would tell us that this term would otherwise
be $D_1(E/k)$ which vanishes as $E$ is smooth. Since $E$ is a smooth
divisor on a smooth variety $\wt X$, $\sL =  {\frac {\sI} {{\sI}^2}}$ is
a line bundle on $E$ and ${\Omega}^1_{E/k}$ is clearly a vector bundle on
$E$. This implies that ${\ov {{\Omega}^1_{nE/k}}}$ is also a vector bundle
on $E$ and for $i \ge 1$, the obvious filtration of 
${\wedge}^i_E\left({\ov {{\Omega}^1_{nE/k}}}\right) =
{\ov {{\Omega}^i_{nE/k}}}$ in terms of the tensor product of the exterior
powers of $\sL$ and ${\Omega}^1_{E/k}$ gives us an exact sequence
\begin{equation}\label{eqn:prelim05}
0 \to {\Omega}^{i-1}_{E/k} {\otimes}_E {\sL} \to
{\ov {{\Omega}^i_{nE/k}}} \to {\Omega}^i_{E/k} \to 0.
\end{equation}
Tensoring this with ${\sL}^j$ and observing that ${\sL}^j = 
{\frac {{\sI}^j} {{\sI}^{j+1}}}$ for $j \ge 1$, we get exact sequence
\begin{equation}\label{eqn:prelim051}
0 \to {\Omega}^{i-1}_{E/k} {\otimes}_E {\sL}^{j+1} \to
{\ov {{\Omega}^i_{nE/k}}} {\otimes}_E {\sL}^{j} \to {\Omega}^i_{E/k} 
{\otimes}_E {\sL}^{j} \to 0.
\end{equation}
Using the right exactness of $H^d$ as before and our assumption, we
get ~\ref{eqn:prelim03} and hence ~\ref{eqn:prelim01}.

Finally, to prove the lemma, we consider the commutative diagram
of exact sequences 
\[
\xymatrix{
0 \ar[r] & {\Omega}^i_{(nE, E)/k} \ar[r] \ar[d] & 
{\Omega}^i_{{nE}/k} \ar@{->>}[d] \ar[r] &
{\Omega}^i_{E/k} \ar[r] \ar@{=}[d] & 0 \\
0 \ar[r] & {\Omega}^{i-1}_{E/k} {\otimes}_E {\sL} \ar[r] & 
{\ov {{\Omega}^i_{nE/k}}} \ar[r] & {\Omega}^i_{E/k} \ar[r] & 0,}
\]
where the bottom row is the exact sequence of ~\ref{eqn:prelim05}.
Observing that the kernel of the middle vertical arrow is a quotient of
${\Omega}^i_{{nE}/k} {\otimes} {\frac {\sI}{{\sI}^n}}$,
a diagram chase above gives an exact sequence
\[{\Omega}^i_{{nE}/k} {\otimes} {\frac {\sI}{{\sI}^n}} \to
{\Omega}^i_{(nE, E)/k} \to {\Omega}^{i-1}_{E/k} {\otimes}_E {\sL} \to 0. \]
Now the lemma follows from the right exactness of $H^d$ on $nE$, the above
claim  ~\ref{eqn:prelim01} and our assumption.
\end{proof}
\begin{lem}\label{lem:prelim2}
Under the conditions of Lemma~\ref{lem:prelim0}, the natural map
\[H^d\left(nE, {\Omega}^i_{{\wt X}/k} {\otimes}_{\wt X} {\sO}_{nE}\right)
\to  H^d\left(E, {\Omega}^i_{E/k}\right)\]
is an isomorphism for $i \ge 0$ and $n \ge 1$.
\end{lem}
\begin{proof} The proof follows from the arguments very similar to that
in the proof of Lemma~\ref{lem:prelim0}. We give a sketch.
First we show this when $n = 1$. In this case, the argument used in the 
proof of the exact sequence ~\ref{eqn:prelim05} also gives the
exact sequence
\[0 \to {\Omega}^{i-1}_{E/k} {\otimes}_E {\sL} \to
{\Omega}^i_{{\wt X}/k} {\otimes}_{\wt X} {\sO}_E 
\to {\Omega}^i_{E/k} \to 0.\]
Taking the exact sequence of $H^d$ on $E$, our hypothesis now proves the
case $n = 1$.

For $n \ge 2$, we can use the $n = 1$ case to reduce the proof of the
lemma to showing that for $n \ge 2$ the natural map
\begin{equation}\label{eqn:prelim10}
H^d\left(nE, {\Omega}^i_{{\wt X}/k} {\otimes}_{\wt X} {\sO}_{nE}\right)
\to H^d\left(E, {\Omega}^i_{{\wt X}/k} {\otimes}_{\wt X} {\sO}_{E}\right)
\end{equation}
is an isomorphism.
We have the short exact sequence of sheaves on $\wt X$
\[0 \to {\sI}/{{\sI}^n} \to {\sO}_{nE} \to {\sO}_E \to 0,\]
which in turn gives the exact sequence
\[0 \to {\Omega}^i_{{\wt X}/k} {\otimes}_{\wt X}  {\sI}/{{\sI}^n} \to
{\Omega}^i_{{\wt X}/k} {\otimes}_{\wt X} {\sO}_{nE} \to
{\Omega}^i_{{\wt X}/k} {\otimes}_{\wt X} {\sO}_{E} \to 0,\]
as ${\Omega}^i_{{\wt X}/k}$ is a vector bundle on $\wt X$.
Taking the exact sequence of $H^d$ on $nE$, we only need to show that
$H^d\left(nE, {\Omega}^i_{{\wt X}/k} {\otimes}_{\wt X} 
{\sI}/{{\sI}^n}\right) = 0$  for $n \ge 2$. But this is proved exactly in the
same way as we proved ~\ref{eqn:prelim01}.
\end{proof}
\begin{lem}\label{lem:prelim1}
Under the conditions of Lemma~\ref{lem:prelim0}, one has
\[H^d\left(nE, {\Omega}^i_{(nE, E)/{\Q}}\right) = 0 
 \ {\rm for} \ i \ge 0 \ {\rm and} \ n \ge 1.\]
\end{lem}
\begin{proof} The spectral sequence of ~\ref{eqn:generalfield} gives a
finite filtration 
\[
{\Omega}^i_{{nE}/{\Q}} = F^0_n \supseteq F^1_n \supseteq \cdots \supseteq
F^{i}_n \supseteq F^{i+1}_n = 0\]
of ${\Omega}^i_{{nE}/{\Q}}$ and a morphism of filtered modules
${\Omega}^i_{{nE}/{\Q}} \surj {\Omega}^i_{{(n-1)E}/{\Q}}$ such that
${\Omega}^i_{{k}/{\Q}} {\otimes}_k {\Omega}^{i-j}_{{nE}/k} \surj
{F^j_n}/{F^{j+1}_n}$ for $0 \le j \le i$.  We also see immediately from the
spectral sequence ~\ref{eqn:generalfield} and by a descending induction on 
$j$ that $F^j_n \surj F^j_1$ for $0 \le j \le i$ and $n \ge 1$. 
In particular, ${\Omega}^i_{(nE, E)/{\Q}}
= {\rm Ker}\left({\Omega}^i_{{nE}/{\Q}} \surj {\Omega}^i_{E/{\Q}}\right)$
has a filtration 
\[{\Omega}^i_{(nE, E)/{\Q}} = F^0_{(n,1)} \supseteq F^1_{(n,1)}
 \supseteq \cdots F^i_{(n,1)} \supseteq F^{i+1}_{(n,1)} = 0\]
with $F^j_{(n,1)} = {\rm Ker}\left(F^j_n \surj F^j_1\right)$
and 
\[{F^j_{(n,1)}}/{F^{j+1}_{(n,1)}} = {\rm Ker}\left({F^j_n}/{F^{j+1}_n}
\surj {F^j_1}/{F^{j+1}_1}\right) \ {\rm for} \ 0 \le j \le i.\]
On the other hand, the spectral sequence  ~\ref{eqn:generalfield} also
shows that ${\Omega}^i_{{k}/{\Q}} {\otimes}_k {\Omega}^{i-j}_{{E}/k} \cong
{F^j_1}/{F^{j+1}_1}$ for $0 \le j \le i$ as $E$ is smooth, and this
gives 
\[{\Omega}^i_{{k}/{\Q}} {\otimes}_k {\Omega}^{i-j}_{{(nE,E)}/k} \surj
{F^j_{(n,1)}}/{F^{j+1}_{(n,1)}}.\]
This in turn gives a surjection
\begin{equation}\label{eqn:prelim10}     
{\Omega}^i_{{k}/{\Q}} {\otimes}_k 
H^d\left(nE, {\Omega}^{i-j}_{(nE, E)/{k}}\right) \surj
H^d\left(nE, {F^j_{(n,1)}}/{F^{j+1}_{(n,1)}}\right)
\ {\rm for} \ 0 \le j \le i \ {\rm and} \ n \ge 1.
\end{equation}
Now the lemma follows by using the exact sequence
\[H^d\left(nE, {F^{j+1}_{(n,1)}}\right) \to 
H^d\left(nE, {F^j_{(n,1)}}\right) \to 
H^d\left(nE, {F^j_{(n,1)}}/{F^{j+1}_{(n,1)}}\right) \to 0,\]
~\ref{eqn:prelim10}, Lemma~\ref{lem:prelim0}
and a descending induction on $j$.
\end{proof}
{\bf Proof of Theorem~\ref{thm:cone}:}
We first show at once that $(i)$ implies $(iv)$ and $(v)$.
We follow the diagram ~\ref{eqn:blowup} and the notations as in the
beginning of this section.
Since $\wt X$ (resp $Z$) is an ${\A}^1$-bundle (resp a ${\P}^1$-bundle)
over $Y$, it is easy to show that $CH^{d+1}(\wt X) = 0$ and
$CH^{d+1}(Z) \cong CH^d(Y)$. Thus it suffices to show that
the natural maps 
\begin{equation}\label{eqn:cone0}
CH^{d+1}(X) \to CH^{d+1}(\wt X) \ {\rm and} \
CH^{d+1}(\ov X) \to CH^{d+1}(Z)
\end{equation}
are isomorphisms.   
We prove the first isomorphism. The proof of the second isomorphism
is exactly the same, once we observe that $E \inj {\wt X} \inj Z$ is the 
exceptional divisor for both $p$ and $\ov p$. We only need to verify
the two conditions of Corollary~\ref{cor:sufficient}.
The first condition is obvious since the inclusion $E \inj {\wt X}$
has a section.
To prove the second condition, it is enough to show that
$H^d\left(nE, {\Omega}^i_{(nE, E)/{\Q}}\right) = 0,$
since dim$(nE) = d$. By Lemma~\ref{lem:prelim1}, this further reduces
to showing that
\begin{equation}\label{eqn:cone1}
H^d\left(E, {\Omega}^i_{E/k} {\otimes}_E {\frac {{\sI}^n}
{{\sI}^{n+1}}}\right)
= 0 \ {\rm for} \ i \ge 0 \ {\rm and} \ n \ge 1.
\end{equation}

Let $T \subset Y$ be a hyperplane section of $Y$ for the given
embedding $Y \inj {\P}^N_k$. Then for $n \ge 1$, one has as short exact 
sequence
\[0 \to {{\sO}_Y}(n-1) \to {{\sO}_Y}(n) \to
{{\sO}_T}(n) \to 0.\] 
Since $T$ is $(d-1)$-dimensional, our assumption now immediately implies 
that $H^d(Y, {{\sO}_Y}(n)) = 0 \ {\forall} \ n \ge 1$.
However, since $\pi$ is an ${\A}^1$-bundle over $Y$
associated to the line bundle ${{\sO}_Y}(1)$ with a section $E$, 
we have ${{\sI}^{n}}/{{\sI}^{{(n+1)}}} \cong {{\sO}_Y}(n)$
for all $n \ge 1$. Hence we get
\begin{equation}\label{eqn:E}
H^d\left(E, {{\sI}^n}/{{\sI}^{{(n+1)}}}\right) = 0 
\ {\forall} \ n \ge 1.
\end{equation} 
This proves ~\ref{eqn:cone1} for $i = 0$. For $i \ge 1$, 
we have $H^d\left(E, {\Omega}^i_{E/k} {\otimes}_E 
{{\sI}^n}/{{\sI}^{{(n+1)}}}\right) \cong
H^d\left(Y, {\Omega}^i_{Y/k}(n)\right)$ under the 
isomorphism $E \xrightarrow {\cong} Y$. But this last group is zero for 
$n \ge 1$ and $i \ge 1$ by the Akizuki-Nakano vanishing theorem 
({\sl cf.} \cite{EV}, Theorem~1.3). This proves ~\ref{eqn:cone1}
and hence $(i)$ implies $(iv)$ and $(v)$.

Now we assume that $k$ is a universal domain. In this case, the
implication $(iv) \Rightarrow (i)$ was shown by Srinivas in
\cite{Srinivas2} (Corollary~2).  Before we prove the other implications,
we first observe that as $\ov X$ has only isolated singularities,
the map $H^{d+1}\left({\ov X}, {\Omega}^i_{{\ov X}/k}\right) \to
H^{d+1}\left({\ov X}, {\ov p}_*\left({\Omega}^i_{Z/k}\right)\right)  
$ is an isomorphism for $i \ge 0$. Now the Leray spectral sequence
gives us for $i \ge 0$, an exact sequence
\begin{equation}\label{eqn:cone2}
\xymatrix@C.5pc{
H^{d}\left(Z, {\Omega}^i_{Z/k}\right) \ar[r] &
~{\underset {n}{\varprojlim}}~H^d\left(nE, {\Omega}^i_{Z/k} {\otimes}_Z
{\sO}_{nE}\right) \ar[r] &  
H^{d+1}\left({\ov X}, {\Omega}^i_{{\ov X}/k}\right) \ar[r] &
H^{d+1}\left(Z, {\Omega}^i_{Z/k}\right) \ar[r] & 0.}
\end{equation}
Furthermore, since ${\ov {\pi}}$ is a ${\P}^1$-bundle, we can use the
Leray spectral sequence again to get
\begin{equation}\label{eqn:cone3}
H^i\left(Z, {\sO}_Z\right) \cong H^i\left(Y, {\sO}_Y\right) \cong 
H^i\left(E, {\sO}_E\right) \ {\forall} \ i \ge 0 \ {\rm and}
\end{equation}
\[H^i\left(Z, {\Omega}^r_{Z/k}\right) \surj  
H^i\left(E, {\Omega}^r_{E/k}\right) \ {\forall} \ i \ge 0 \
r \ge 1.\]

{\sl Proof of} $(i) \Leftrightarrow (ii)$. 
Suppose $H^{d+1}\left({\ov X}, {\sO}_{\ov X}\right) = 0$. Then $(i)$
follows from ~\ref{eqn:cone2} for $i = 0$ and ~\ref{eqn:cone3},
once we observe that $E \inj nE$ has a section for all $n \ge 1$.
If $(i)$ holds, then we have already seen in ~\ref{eqn:prelim01} that
$H^d\left(nE, {\sI}/{{\sI}^n}\right) = 0$ for $n \ge 1$ and hence
$H^d\left(nE, {\sO}_{nE}\right) \cong H^d\left(E, {\sO}_{E}\right)$
for all $n \ge 1$. Now $(ii)$ follows from  
~\ref{eqn:cone2} for $i = 0$ and ~\ref{eqn:cone3}.

{\sl Proof of} $(ii) \Leftrightarrow (iii)$. We only need to show that 
$(ii) \Rightarrow (iii)$. Since 
$H^{d+1}\left(W, {\Omega}^i_{W/k}\right)$ are birational invariants of
$W$, we can replace $W$ by $Z$ everywhere below. 
It suffices then
to show that $H^{d+1}\left({\ov X}, {\Omega}^i_{{\ov X}/k}\right) \cong
H^{d+1}\left(Z, {\Omega}^i_{Z/k}\right)$ for $i \ge 1$.
By ~\ref{eqn:cone2} and ~\ref{eqn:cone3}, this reduces to showing that 
$H^d\left(nE, {\Omega}^i_{{nE}/k} {\otimes}_{nE} {\sI}/{{\sI}^n}\right)
= 0$ for $n \ge 1$. But this follows from the Akizuki-Nakano vanishing 
theorem and the claim ~\ref{eqn:prelim01} in the proof of 
Lemma~\ref{lem:prelim0}.

{\sl Proof of} $(iii) \Leftrightarrow (iv)$.
We have already shown that $(iii) \Leftrightarrow (i)$ and 
$(i) \Leftrightarrow (iv)$. So we get $(iii) \Leftrightarrow (iv)$.

{\sl Proof of} $(iv) \Leftrightarrow (v)$.
We have seen before that $(iv) \Leftrightarrow (i) \Rightarrow (v)$.
So we are only left to show that 
$(v) \Rightarrow (iv)$. \\
For this, it suffices to show that ${\rm Ker}\left(CH^{d+1}(\ov X)
\to CH^{d+1}(Z)\right) \surj CH^{d+1}(X)$.
Let $H \inj {\ov X}$ be the hyperplane at infinity, i.e., $H$ is the
complement of $X$ in $\ov X$ and is isomorphic to $Y$. 
Moreover, the inclusion $H \inj {\ov X}$ factors through the inclusion
$H \xrightarrow {j} Z$ and the natural map $CH^d(H) \xrightarrow {j_*}
CH^{d+1}(Z)$ is an isomorphism. Now we have the following commutative
diagram
\[
\xymatrix@C.6pc{
& & F \ar@{->>}[r] \ar[d] & CH^d(H) \ar[d]^{j_*} & \\
0 \ar[r] & {\rm Ker}\left({\ov p}^*\right) \ar[d] \ar[r] &
CH^{d+1}(\ov X) \ar@{->>}[d] \ar[r] & CH^{d+1}(Z) \ar[r] \ar[d] & 0 \\
0 \ar[r] & {\rm Ker}\left({p}^*\right) \ar[r] &
CH^{d+1}(X) \ar[r] & CH^{d+1}(\wt X) \ar[r] & 0, }
\]  
where $F$ is the free abelian group on the closed points of $H$.
Note that as $H$ is contained in the smooth locus of $\ov X$,
the map $F \to H^{d+1}(\ov X)$ is well defined and the composite
$F \to H^{d+1}(\ov X) \to CH^{d+1}(X)$ is clearly zero.
We already know that $CH^{d+1}(\wt X) = 0$. Since $j_*$ is an isomorphism,
a diagram chase shows that ${\rm Ker}\left({\ov p}^*\right) \surj
CH^{d+1}(X)$. This completes the proof of Theorem~\ref{thm:cone}. 
$\hfill \square$

School of Mathematics \\
Tata Institute Of Fundamental Research \\
Homi Bhabha Road \\
Mumbai,400005, India \\
email : amal@math.tifr.res.in \\
\end{document}